\documentclass[14]{article}
\usepackage{arxiv}
\usepackage[utf8]{inputenc} % allow utf-8 input
\usepackage[T1]{fontenc}    % use 8-bit T1 fonts
\usepackage{amsmath}
\usepackage{hyperref}       % hyperlinks
\usepackage{url}            % simple URL typesetting
\usepackage{booktabs}       % professional-quality tables
\usepackage{amsfonts}       % blackboard math symbols
\usepackage{nicefrac}       % compact symbols for 1/2, etc.
\usepackage{microtype}      % microtypography
\usepackage{lipsum}
\usepackage{cases}
 \usepackage{setspace}
\usepackage{graphicx}   
\usepackage{subcaption} 
\usepackage{caption}
\usepackage{tikz, pgfplots}
\usepackage{tikz}
\usetikzlibrary{3d}
\usetikzlibrary{hobby}
\usepackage{amssymb}
\usepackage{mathrsfs}
\usepackage{latexsym}
\usepackage{enumerate}
\usepackage{color}
\usepackage{ulem}
\usepackage{amsfonts}
\usepackage{verbatim}
\usepackage{amsthm}
\usepackage{dsfont}
\usepackage{url}
\usepackage{tocloft}
\usepackage{youngtab} % for Y-diagrams
 \usepackage{ytableau} % for Y-diagrams
\usepackage{enumerate}
\usepackage{mathtools}
\usepackage{multirow}
\usepackage{array}
\usepackage{xpatch}

\usepackage{tabularx}
\usepackage{makecell}

\allowdisplaybreaks
\usepackage{cellspace} % Load the cellspace package
\usepackage{slashbox}
\usepackage{diagbox}
\pgfplotsset{compat=1.18}

\xapptocmd{\proof}{\mbox{}\par\nobreak}{}{}
\newtheorem{theorem}{Theorem}[section]

\newtheorem{condition}[theorem]{Condition}

\newtheorem{prop}[theorem]{Proposition}

\newtheorem*{authorcontributions}{Author Contributions}
\newtheorem*{conflictofinterest}{Conflict of Interest}

\newtheorem*{funding}{Funding}
\numberwithin{equation}{section}

\newcolumntype{C}[1]{>{\centering\arraybackslash}p{#1}}

\title{A Point on Discrete vs Continuous State-Space Markov Chains}

\author{
  Mathias Muia\thanks{Corresponding author: \texttt{mnmuia@southalabama.edu}} \\
  Department of Mathematics and Statistics \\
  University of South Alabama, Mobile, AL 36688, USA
  \and
  Martial Longla\thanks{\texttt{mlongla@olemiss.edu}} \\
  Department of Mathematics \\
  University of Mississippi, University, MS 38677, USA
}

\date{}

\renewcommand{\arraystretch}{1.5}  % Increase the row height by 1.5 times
\begin{document}

\maketitle

\begin{abstract}
This paper investigates the effects of discrete marginal distributions on copula-based Markov chains. We establish results on mixing properties and parameter estimation for a copula-based Markov chain model with Bernoulli($p$) marginal distributions, emphasizing some distinctions between continuous and discrete state-space Markov chains. We derive parameter estimators utilizing the maximum likelihood estimation (MLE) method and explore alternative estimators of $p$ that are asymptotically equivalent to the MLE. Furthermore, we provide the asymptotic distributions of these parameter estimators. A simulation study is conducted to evaluate the performance of the various estimators for $p$. Additionally, we employ the likelihood ratio test to assess independence within the sequence.

\end{abstract}

\section{Introduction}
Markov chain models have been widely used in the literature to represent the relationship between consecutive observations in a sequence of Bernoulli trials. Johnson and Klotz (1974) \cite{johnson1974atom} employed a Markov chain generalization of a binomial model to analyze the crystal structure of a two-metal alloy in their metallurgical studies. Crow (1979) \cite{crow1979approximate} investigated the applications of two-state Markov chains in telecommunications in his study focusing on approximating confidence intervals. Brainerd and Chang (1982) \cite{brainerd1982number} utilized two-state Markov chains to address problems in linguistic analysis. These examples illustrate that Markov chain models can provide a solid foundation for modeling a wide range of two-state sequential phenomena.

Moreover, in the study of sequences consisting of different alternatives, there has been notable interest in using group or run tests to evaluate the randomness of the sequences. These tests are based on the count of changes between outcomes in the sequence, aiming to assess whether the sequence displays randomness (denoted as the null hypothesis, $H_0$) or exhibits a type of dependence similar to simple Markov chains (alternative hypothesis, $H_1$). These counts, denoted by $n_{ij}\ (\text{where} \ i,j=0,\cdots, s\ \text{denotes the states taken by the sequence})$ are commonly referred to as \textit{transition numbers} of the sequence. The transition numbers have been widely used in literature to construct run tests and in estimation of model parameters (see David (1947b) \cite{david1947power}, Goodman (1958) \cite{goodman1958simplified}, Anderson and Goodman (1957) \cite{anderson1957statistical} and Billingsley (1961a) \cite{billingsley1961statistical} for details). As used in this work, transition numbers prove to be very useful in computing the maximum likelihood estimators of parameters. 

\subsection{Copulas and Markov chains}
A 2-copula is a bivariate function $C: [0,1]^{2}\rightarrow [0,1]$ such that $C(0,u)=C(u,0)=0$, $C(u,1)=C(1,u)=u, \forall \ u\in [0,1]$ and $C(u_{1},v_{1})+C(u_{2},v_{2})-C(u_{1},v_{2})-C(u_{2},v_{1})\geq 0, \forall\ [u_{1},u_{2}]\times[v_{1},v_{2}]\subset [0,1]^{2}$. For some $k\geq 2,$ let $C_{1}, \cdots, C_{k}$ be a set of $k$ copulas and $a_{1}, \cdots, a_{k}$ a set of real constants such that $a_{1}+\cdots+a_{k}=1$ and $0\leq a_{i}\leq 1$ for each $i=1,\cdots, k$, then the sum given by 
$ C(u,v)=a_1C_1(u,v)+\cdots+a_kC_k(u,v) $
is a copula referred to as the convex combination of $C_{1}, \cdots, C_{k}$ (Longla, 2015) \cite{longla2015mixtures}.

The copula $C(F(x),G(y))=H(x,y)$, of the random variables $X$ and $Y$, is known to be unique if $F$ and $G$ are continuous; otherwise $C$ is uniquely determined on $Range(F)\times Range(G)$ (Nelsen (2006) \cite{nelsen2006methods}). In the recent past, copulas have gained popularity in the modeling of temporal
dependence within Markov processes. This is because, when used to represent joint distributions of random variables, they enable the detection of possible correlations or links between variables (Darsow, Nguyen and Olsen (1992) \cite{darsow1992copulas}). A copula-based Markov chain is a stationary Markov chain generated by the copula of its consecutive states and an invariant distribution.  Copula-based Markov chains enable the exploration of dependence and mixing properties inherent in Markov processes.   This makes them an invaluable tool for investigations and allows for the study of scale-free measures of dependence while remaining invariant under monotonic transformations of the variables. Consequently, copula methods have drawn a lot of interest across diverse fields, including technology, finance, and the natural world.

This paper focuses on a stationary Markov chain based on a copula from the Frechet family of copulas. The Fr\'echet family of copulas has the form \begin{align}\label{Frechet}
		C(u,v)=aM(u,v)+(1-a-b)\Pi(u,v)+bW(u,v), \quad a,b\geq 0,\ 0\leq a+b\leq 1,
\end{align}
where $\Pi(u,v)=uv$ is the independence copula and  
$M(u,v)=\min(u,v)$ and $W(u,v)=\max(u+v-1,0)$ are the Fr\'echet-Hoeffding upper and lower bounds respectively. We study two cases: a continuous state-space Markov chain with uniform marginal distributions and a discrete state-space Markov chain with Bernoulli($p$)  marginal distributions. Maximum likelihood estimators (MLEs) of model parameters are derived in both cases and their asymptotic distributions are provided. The likelihood ratio test is used to test for randomness in the model. 

Studies involving inferential statistics for discrete cases closely similar to the one considered here have been presented in the works of Billingsley (1961a) \cite{billingsley1961statistical}, Billingsley (1961b) \cite{billingsley1961methods}, Goodman (1958) \cite{goodman1958simplified}, Klotz (1972) \cite{klotz1972markov}, Klotz (1973) \cite{klotz1973statistical}, Price (1976) \cite{price1976note}, Lindqvist (1978) \cite{lindqvist1978note},  and others. Goodman (1958) \cite{goodman1958simplified} considered a single sequence of alternatives consisting of a long chain of observations to derive some long sequence group tests. Klotz (1973) \cite{klotz1973statistical} presented small and large sample distribution theory for the sufficient statistics of a model for Bernoulli trials with Markov dependence. He provided estimators for the model parameters and showed that the uniform most powerful (u.m.p.) unbiased test of independence is the run test. On the other hand, Price (1976) \cite{price1976note} and Lindqvist (1978) \cite{lindqvist1978note} focused on estimating the model parameters only. Price (1976) \cite{price1976note} used Monte Carlo techniques to investigate finite sample properties for the parameters of a dependent Bernoulli process.   Lindqvist (1978) \cite{lindqvist1978note} considered the weaker assumption of non-Markovian dependence and showed that the MLE, $(\hat{p},\hat{a})$ is a strongly consistent, and asymptotically normally distributed estimator for $(p, a)$. In the context of Klotz (1973) and Lindqvist (1978), as well as in this paper, $p$ represents the Bernoulli (frequency) parameter of a sequence of Bernoulli trials, while $a$ is an additional dependence parameter. Lindqvist (1978) \cite{lindqvist1978note} further showed that $(\hat{p},\hat{a})$ is asymptotically equivalent to $(\Bar{p},\bar{a})$, where $\Bar{p}$ is the sample mean and $\bar{a}$ is the empirical correlation coefficient of $X_{t-1}$ and $X_t$.

The rest of the paper is organized as follows: In Section \ref{begin here}, we introduce our model and detail the transition probabilities of the discrete Markov chain model considered. Mixing properties of the Markov chain model are discussed in Subsection \ref{Mixing coefficients}. Section \ref{Parameter estimation and CLTs} covers parameter estimation, where we derive MLEs and their asymptotic distributions. In Section \ref{test for independence} we provide a test for independence in the sequence. The paper concludes with a simulation study in Section \ref{Simulations Study}.

\par
\section{The model}\label{begin here}
Consider a stationary Markov chain model based on the copula 
\begin{align}\label{our model}
	C(u,v)=a M(u,v)+(1-a)W(u,v)   ,   \qquad a\in(0,1).       \end{align}
The copula (\ref{our model}) is from the Frechet/Mardia family of copulas given by (\ref{Frechet}) with $a+b=1$. Note that a Mardia copula with a parameter $\theta\in[-1,1]$ is a Frechet copula with $a + b = \theta^2$. For details on these copula families, see Nelsen (2006) \cite{nelsen2006methods}, Joe (1997) \cite{joe1997multivariate}, or Durante and Sempi (2015) \cite{durante2015principles}. 

Longla (2014) \cite{longla2014dependence} showed that a stationary Markov chain with copula (\ref{our model}) and uniform marginals has $n$-steps joint cumulative distribution function (cdf) of the form \begin{align}\label{uniform01}
	C_n(u,v)=    1/2[1+(2a-1)^n]M(u,v)+1/2[1-(2a-1)^n]W(u,v). 
\end{align}

If the marginals are  Bernoulli($p$), the joint distributions can be characterized by the transition probabilities of the Markov chain and the stationary distribution. Note that Sklar’s theorem states that C is only identifiable on $ Ran(F) \times Ran (G).$ Theorem \ref{transition matrices} below provides the transition matrices for the Markov chain. 

\begin{theorem}\label{transition matrices}
	Let $X=(X_t)_{t=0,\cdots, n}$ be a stationary Markov chain generated by copula (\ref{our model})
	and Bernoulli($p$) marginals. Then $X$ has transition matrices given by
	\begin{equation}\label{tmatrples}
		A =  \begin{pmatrix}
			\frac{ap+1-2p}{1-p} & \frac{p(1-a)}{1-p}\\
			1-a & a
		\end{pmatrix}, \  p<1/2
	\end{equation}
	and
	\begin{equation}\label{tmatrpgeq}
		A =  \begin{pmatrix}
			a & 1-a\\
			\frac{(1-p)(1-a)}{p} & \frac{a(1-p)+2p-1}{p}\\
		\end{pmatrix}, \ p\geq 1/2.
	\end{equation}
\end{theorem}  
\begin{proof}
	Due to stationarity of the Markov chain, $P(X_t=1)=1-P(X_t=0)=p$. It suffices to determine $P(X_0=0, X_1=0)$. \begin{align*}
		P(X_0=0,X_1=0)&=P(X_0\leq 0,X_1\leq 0)=C(1-p,1-p)\\ &=a(1-p)+(1-a)\max(1-2p,0)\\\nonumber
		&=\begin{cases}
			ap+1-2p,&\quad\quad \text{if}\ p<1/2\\
			a(1-p),&\quad\quad\text{if}\ p\geq 1/2,
		\end{cases}
	\end{align*}  where the second equality follows from Sklar's Theorem (Sklar (1959)) \cite{sklar1959fonctions}. The rest of the probabilities are computed using the hypothesis that the Markov chain is stationary with Bern(p) marginals. For $p<1/2$, the joint distribution is given by Table \ref{table1} below. 
\begin{table}[htpb]
\centering
\def\arraystretch{1.0}
\begin{tabular}{|c||c|c|c|}
\hline
& \multicolumn{2}{c|}{$X_1$} & \\ 
\cline{2-3}
$X_0$ & $0$ & $1$ & $P_{X_{0}}(x)$ \\
\hline\hline
$0$ & $ap + 1 - 2p$ & $p(1 - a)$ & $1 - p$ \\
\hline
$1$ & $p(1 - a)$ & $ap$ & $p$ \\
\hline
$P_{X_{1}}(y)$ & $1 - p$ & $p$ & $1$ \\
\hline
\end{tabular}
\caption{Distribution of $(X_0, X_1)$ when $p < 1/2$.}
\label{table1}
\end{table}

Using the definition of conditional probability, $P(X_{1}=j|X_{0}=i)=P(X_{1}=j,X_{0}=i)/P(X_{0}=i)$, in Table \ref{table1} we obtain the transition matrix (\ref{tmatrples}). A similar argument is used for $p\geq 1/2$ to get the transition matrix (\ref{tmatrpgeq}).
\end{proof}

The statistical behavior of the distribution at \( p = 1/2 \) reflects a symmetric transition structure where the probability of transitioning between states is balanced by the parameter \( a \). This symmetry simplifies the interpretation and analysis of the Markov chain, making it an interesting case for studying the dynamics governed by the mixture copula model.

The transition matrix (\ref{tmatrples}) is similar to the one considered by Klotz (1973) \cite{klotz1973statistical}. In our case, we have a frequency parameter $p=P(X_t=1)$ and an additional dependence parameter $a=P(X_t=1|X_{t-1}=1)$, when $p<1/2$ or $a=P(X_t=0|X_{t-1}=0)$, when $p\geq 1/2$. In Klotz (1973) \cite{klotz1973statistical}, the condition $\max(0,(2p-1)/p)\leq a\leq 1$ was required so that the transition probabilities of the model remain bounded between $0$ and $1$. The model presented in this work does not require this condition because the copula conditions sufficiently define the transition probabilities based on whether $p<1/2$ or $p\geq 1/2$. Specifically, the transition probabilities remain bounded between 0 and 1 for any value of $0<a<1$, and they satisfy $\sum\limits_j p_{ij}=1$ for $j=0,1$.

\begin{prop}\label{nth transition matrices}
	The $n^{th}$ transition matrices for (\ref{tmatrples}) and (\ref{tmatrpgeq}) are respectively\begin{align}\label{nthtranp<1/2}
		A^{n} =  \begin{pmatrix}
			1-p+p(\frac{a-p}{1-p})^{n} & p-p(\frac{a-p}{1-p})^{n}\\
			1-p-(1-p)(\frac{a-p}{1-p})^{n} & p+(1-p)(\frac{a-p}{1-p})^{n}
		\end{pmatrix}, \quad \quad p<1/2
	\end{align}
	and
	\begin{align}\label{nthtranp>1/2}
		A^{n} =  \begin{pmatrix}
			1-p+p\left(\frac{a+p-1}{p}\right)^n & p-p\left(\frac{a+p-1}{p}\right)^n\\
			1-p-(1-p)\left(\frac{a+p-1}{p}\right)^n & p+(1-p)\left(\frac{a+p-1}{p}\right)^n
		\end{pmatrix},\quad \quad p\geq 1/2
	\end{align}
\end{prop}
\begin{proof}
	For $p<1/2$, the eigenvalues of $A$ are obtained by solving the characteristic equation:
	\begin{equation*}
		\lambda^{2}-\frac{a+1-2p}{1-p}\lambda+\frac{a-p}{1-p}=0.
	\end{equation*}
	Solving this equation yields $\lambda_{1}=1$ and $\lambda_{2}=\frac{a-p}{1-p}$. Then,
	\begin{equation}\label{Ales}
		A^{n} =  \begin{pmatrix}
			a_{00}+b_{00}\lambda_{2}^{n} & a_{01}+b_{01}\lambda_{2}^{n}\\
			a_{10}+b_{10}\lambda_{2}^{n} & a_{11}+b_{11}\lambda_{2}^{n}
		\end{pmatrix},
	\end{equation}
	where $a_{ij}$, $b_{ij}$, $i,j=0,1$ are coefficients to be determined, and $\lambda_{2}$ is the non-simple eigenvalue of $A$.
	By using (\ref{Ales}) to write $A^{0}$ and $A^{1}$, we can solve for $a_{ij}$ and $b_{ij}$, $i,j=0,1$, by forming four simultaneous equations, which are easy to solve. This yields the $n^{th}$ transition matrix (\ref{nthtranp<1/2}).
	Similar steps are used to obtain the $n^{th}$ transition matrix for (\ref{tmatrpgeq}).
\end{proof}

\subsection{Mixing}\label{Mixing coefficients}  
Most of the dependence and mixing coefficients are heavily influenced by the copula of the model of interest (Longla and Peligrad (2012) \cite{longla2012some}). Longla (2014) \cite{longla2014dependence} proposed a set of copula families that generate exponential $\rho$-mixing Markov chains. He showed that stationary Markov chains generated by copulas from the the Frechet(Mardia) families and uniform marginals are exponentially $\phi$-mixing, therefore $\rho$-mixing and geometrically $\beta$-mixing  if $a+b\neq 1$. However, when $a+b=1$, mixing does not occur (Longla, 2014) \cite{longla2014dependence}. Furthermore, Longla (2015) \cite{longla2015mixtures} showed that copulas from these families fail to generate $\psi$-mixing continuous state-space strictly stationary Markov chains for any $a$ and $b$, achieving at best `lower $\psi$-mixing'. Longla (2015) \cite{longla2015mixtures} also noted the absence of $\psi$-mixing for copulas from Mardia and Frechet families for any $a$ and $b$. Different results emerge when discrete marginal distributions are employed.

We define two mixing coefficients in this paper: the $\phi$-mixing coefficient, introduced by Ibragimov (1959) \cite{ibragimov1959some} and further studied by Cogburn (1960) \cite{cogburn1960asymptotic}, and the $\psi$-mixing coefficient, first introduced by Blum, Hanson, and Koopmans (1963) \cite{blum1963strong} and refined by Philip (1969) \cite{philipp1969central}. In the context of a probability space $(\Omega,\mathcal{F},P)$ with sigma algebras $\mathcal{A},\mathcal{B}\subset\mathcal{F}$, the coefficients are defined in literature as follows (see Bradley (2007) \cite{bradley2007introduction}) \begin{align*}
	\phi(\mathcal{A}, \mathcal{B})=\sup_{A\in\mathcal{A},B\in\mathcal{B},  P(A)>0}\left|\frac{P(B\cap A)}{P(A)}-P(B)\right|,
\end{align*}
\begin{align*}
	\psi(\mathcal{A}, \mathcal{B})=\sup_{A\in\mathcal{A},B\in\mathcal{B},  P(A)>0,P(B)>0}\left|\frac{P(A\cap B)}{P(A)P(B)}-1\right|.
\end{align*}
For random variables $n$-steps apart such that $X_0\in A$ and $X_n\in B$, we define the mixing coefficients 
\begin{equation}    \phi(n)=\sup\limits_{A\in\mathcal{A},B\in\mathcal{B},P(A)>0}\left|\frac{P^n(B\cap A)}{P(A)}-P(B) \right|, \quad 
\end{equation}
\begin{equation}\label{psi}  \psi(n)=\sup\limits_{A\in\mathcal{A},B\in\mathcal{B},P(A)>0,P(B)>0}\left|\frac{P^n(A\cap B)}{P(A)P(B)}-1 \right|, \quad 
\end{equation}
where $P^n(A\cap B)=P(X_0\in A,X_n\in B)$.
A random process is said to be $\psi-$mixing if $\psi(n)\longrightarrow 0\ \text{as}\ n\longrightarrow\infty$   and $\phi-$mixing if $\phi(n)\longrightarrow 0\ \text{as}\ n\longrightarrow\infty$. For the case involving discrete state-space Markov chains generated by the copula given by (\ref{our model}) and Bernoulli$(p)$ marginal distribution, we establish the following theorem:

\begin{theorem}\label{discrete}
	Let $X=(X_t)_{t=0,\cdots, n}$ be a stationary Markov chain generated by copula (\ref{our model})   and Bernoulli($p$) marginal distributions. Then $X$ is exponentially $\psi$-mixing with $\psi(n)
		=\frac{1-p}{p}\left|\frac{a-p}{1-p}\right|^{n},$  
		for $p\geq 1/2$, and $\psi(n)=\frac{p}{1-p}\left|\frac{a+p-1}{p} \right|^n,$ for $p<1/2$. Hence, \( X \) is also exponentially $\phi$-, $\rho$-, $\beta$-, and $\alpha$-mixing.
\end{theorem}
\begin{proof} 
	Consider the coefficient defined by Equation (\ref{psi}) with $\Omega=\{0,1\}$. For $p<1/2$, we examine \begin{align*}\left|\frac{P^n(A\cap B)}{P(A)P(B)}-1\right|,\end{align*} where $A$ and $B$ are chosen only from $\{0\}$ or $\{1\}$ since $P(\{0,1\}\cap A)=P(A)$. The joint distribution of $(X_{0},X_{n})$ when $p<1/2$ is obtained from the $n^{th}$ transition matrix (\ref{nthtranp<1/2}) and using $P(X_{1}=j,X_{0}=i) =P(X_{1}=j|X_{0}=i)P(X_{0}=i)$, where $i,j=0,1$. This is given in Table \ref{table3}.

    \begin{table}[h]
\centering
\def\arraystretch{1.2}
\begin{tabular}{|c||c|c|c|}
\hline
& \multicolumn{2}{c|}{$X_n$} & \\
\cline{2-3}
$X_0$ & $0$ & $1$ & $P_{X_0}(x)$ \\
\hline\hline
0 & $(1 - p)^2 + p(1 - p)\left(\frac{a - p}{1 - p}\right)^n$ 
  & $p(1 - p) - p(1 - p)\left(\frac{a - p}{1 - p}\right)^n$ 
  & $1 - p$ \\
\hline
1 & $p(1 - p) - p(1 - p)\left(\frac{a - p}{1 - p}\right)^n$ 
  & $p^2 + p(1 - p)\left(\frac{a - p}{1 - p}\right)^n$ 
  & $p$ \\
\hline
$P_{X_n}(y)$ & $1 - p$ & $p$ & 1 \\
\hline
\end{tabular}
\caption{Distribution of $(X_0, X_n)$ when $p < 1/2$.}
\label{table3}
\end{table}

	\begin{comment}
	\bgroup
	\def\arraystretch{0.8}
	\begin{table}[h]
		\centering
		\begin{tabular}{|c||*{5}{c|}}\hline
			\backslashbox{$X_{0}$}{$X_{n}$}
			&\makebox[3em]{0}&\makebox[3em]{1}&\makebox[3em]{$P_{X_{0}}(x)$}\\ \hline\hline
			0    & $(1-p)^{2}+p(1-p)(\frac{a-p}{1-p})^{n}$ & $p(1-p)-p(1-p)(\frac{a-p}{1-p})^{n}$ & $1-p$ \\\hline
			1   & $p(1-p)-p(1-p)(\frac{a-p}{1-p})^{n}$ & $p^{2}+p(1-p)(\frac{a-p}{1-p})^{n}$ & p\\\hline
			$P_{X_{n}}(y)$ & $1-p$ & $p$ & 1 \\  \hline
		\end{tabular}
		\caption{Distribution of $(X_{0},X_{n})$ when $p<1/2.$}
		\label{table3}
	\end{table}
	\egroup
\end{comment}
	After computing the four quantities, \begin{align*}
		|\frac{P(X_{0}=i,X_{n}=j)}{P(X_{0}=i)P(X_{n}=j)}-1|
	\end{align*} from Table \ref{table3} we end up with
	\begin{align*}
		\psi(n)&=\sup\lbrace \frac{p}{1-p}\left|\frac{a-p}{1-p}\right|^{n},\left|\frac{a-p}{1-p}\right|^{n},\frac{1-p}{p}\left|\frac{a-p}{1-p}\right|^{n}\rbrace
		=\frac{1-p}{p}\left|\frac{a-p}{1-p}\right|^{n}  
	\end{align*}
	For $p\geq 1/2$, a similar argument leads to $   \psi(n)=\frac{p}{1-p}\left|\frac{a+p-1}{p} \right|^n.
	$
	Therefore, for all $a,p\in (0,1)$, $\psi(n)\longrightarrow 0\ \text{at an exponential rate as}\ n\longrightarrow\infty$. The rest of the proof follows from the fact that $\psi$-mixing implies 
\end{proof}

\section{Parameter estimation asymptotic normality}\label{Parameter estimation and CLTs}
The characteristics of the model parameter estimators, including the shape of their asymptotic distributions, rely on the chosen marginal distribution.
\subsection{Markov chain with uniform marginals}
Assuming $(X_{t})_{t=0,\cdots, n}$ is a Markov chain generated by copula (\ref{our model}) with uniform$(0,1)$ marginal distribution and utilizing Sklar's theorem, the joint cumulative distribution function of $(X_{t}, X_{t+1})$ is represented as:
\begin{equation}\label{joint distribution}
	F(x_{t}, x_{t+1}) = C(F(x_{t}), F(x_{t+1})) = aM(x_{t}, x_{t+1}) + (1-a)W(x_{t}, x_{t+1}).
\end{equation}
The copulas $M$ and $W$ have support on the main diagonal and secondary diagonal of $[0,1]^2$ respectively (Nelsen (2006) \cite{nelsen2006methods}). The joint distribution (\ref{joint distribution}) implies $P(X_{t}=X_{t+1}) = a$ and $P(X_{t}=1-X_{t+1}) = 1-a$ due to the properties of mixture distributions. The indicator random variable $I(X_{t}=X_{t+1}), \ t \in {0,1,\ldots, n-1}$, taking the value $1$ when $X_{t}=X_{t+1}$ and $0$ otherwise, follows a Bernoulli($a$) distribution with probability of success $a$ and probability of failure $1-a$. The parameter $a$ can be estimated by the sample mean:
\begin{align}\label{Yn estimator of a}
	Y_{n} = \frac{1}{n}\sum_{t=0}^{n-1}I(X_t=X_{t+1})\end{align}
as demonstrated in Proposition (\ref{lemma on independence of indicators}) below.

%\begin{comment}
%\newpage
\begin{prop}\label{lemma on independence of indicators}
	Let $(X_t)_{t=0,\cdots,n}$ be a Markov chain generated by copula (\ref{our model}) and the uniform(0,1) distribution. We have the following: \begin{itemize}
		\item[i.] The random variables $I(X_{t}=X_{t+1}), t=0,\cdots, n-1$  are independent and identically distributed (i.i.d) Bernoulli random variables.
		\item[ii.] The random variable (\ref{Yn estimator of a})
		is an unbiased and consistent estimator of $a$. Moreover, $Y_n$ satisfies the Central Limit Theorem.
		\item[iii.]   The MLE of $a$ is the method of moments estimator given by  (\ref{Yn estimator of a}). 
	\end{itemize}  
\end{prop}

\begin{proof}
	\begin{itemize}
		\item[i.]  We prove independence in the sequence by induction. For $m=2$ random variables, let us consider without loss of generality the random variables $I(X_1=X_0)\ \text{and}\ I(X_2=X_1)$. Then 
			\begin{align*}
				P(X_1 = X_0,\ & X_2 = X_1)
				= \int_{0}^{1} P(X_1 = X_0,\ X_2 = X_1 \mid X_0 = x)\ dP_{X_0}(x) \\
				&\quad \text{(by the law of total probability, and since under the copula $M$, the} \\
				&\quad \text{joint distribution of $(X_0, X_1)$ is concentrated on the main diagonal)} \\
				&= \int_{0}^{1} P(X_1 = x \mid X_0 = x)\cdot P(X_2 = x \mid X_0 = x,\ X_1 = x)\ dx \\
				&\quad \text{(by the multiplication rule, and using that \(X_0 \sim \text{Unif}(0,1)\))} \\
				&= \int_{0}^{1} P(X_1 = x \mid X_0 = x)\cdot P(X_2 = x \mid X_1 = x)\ dx \\
				&\quad \text{(by the Markov property)} \\
				&= \int_{0}^{1} a \cdot a\ dx = a^2.
			\end{align*}
			since $P(X_{t+1} = x \mid X_t = x) = a$ for all $x \in [0,1]$ and $t=0,\cdots,n-1$. Therefore,
			\begin{align}\label{indep for pair}
				P(X_1=X_0, X_2=X_1)=P(X_1=X_0)P(X_2=X_1),
			\end{align}
			establishing independence for the base case.
		
		Induction Step: Suppose independence holds for $m=k>2$, that is,   \begin{align}\label{indp for n rvs}
			P(X_{1}=X_{0}, X_{2}=X_{1},\cdots,X_{k+1} = X_{k})= \prod\limits_{t=0}^k P(X_{t+1} = X_{t}) = a^{k+1}.
		\end{align}
		We now prove it for $m=k+1$. Consider  \begin{align*}
			P(X_{1}&=X_{0}, X_{2}=X_{1},\cdots,X_{k+1} = X_{k},X_{k+2} = X_{k+1})\\
			&=P(X_{k+2} = X_{k+1}|X_{k+1} = X_{k}, \cdots, X_{2} = X_{1},X_{1}=X_{0}) \times \\&\qquad P(X_{1}=X_{0}, X_{2}=X_{1},\cdots,X_{k+1} = X_{k})\ \text{(by multiplication rule}\\ &\qquad\text{of probabilities)}\\
			&=P(X_{k+2} = X_{k+1}|X_{k+1} = X_{k})\times \prod\limits_{t=0}^k P(X_{t+1} = X_{t})\ \text{(by the Markov }\\ &\qquad\text{property and (\ref{indp for n rvs}) above)}\\
			&=\prod\limits_{i=0}^{k+1} P(X_{t+1} = X_{t})= a^{k+2}\ \text{(due to (\ref{indep for pair})).}
		\end{align*}
		Thus, by induction, independence holds for all $m$.\\
		
		To show that the random variables \( I(X_{t+1} = X_t) \) are \(\text{Bernoulli}(a)\), note that by the model assumption, for each \( t \), we have \( P(X_{t+1} = X_t) = a \). Therefore, the indicator variable
		\[
		I(X_{t+1} = X_t) =
		\begin{cases}
			1 & \text{if } X_{t+1} = X_t, \\
			0 & \text{otherwise},
		\end{cases}
		\]
		satisfies
		\[
		P(I(X_{t+1} = X_t) = 1) = a, \quad P(I(X_{t+1} = X_t) = 0) = 1 - a,
		\]
		which implies that \( I(X_{t+1} = X_t) \sim \text{Bernoulli}(a) \).
		\item[ii.] The proof follows from the properties of the estimator of the Bernoulli parameter
		for an i.i.d. Bernoulli$(a)$ sample.

		\item[iii.]  Let $Y_t=I(X_{t}=X_{t+1})$, 
		\begin{align}
			L(X_t,X_{t+1})=a^{Y_t}(1-a)^{1-Y_t}.\end{align} The likelihood function is \begin{align}\label{density of xi=xi+1 as a bern(a)}  L(a)=\prod\limits_{t=0}^{n-1}L(X_t,X_{t+1}) =  a^{\sum_{t=0}^{n-1}Y_t}(1-a)^{n-\sum_{t=0}^{n-1}Y_t}.\end{align} 
		It follows that $\hat{a}=\frac{1}{n}\sum_{t=0}^{n-1}Y_t=\frac{1}{n}\sum_{t=0}^{n-1}I(X_t=X_{t+1})$. 
	\end{itemize}
\end{proof}
%\end{comment}

%\newpage
\subsection{\texorpdfstring{Stationary Markov chain with Bern$(p)$ marginal distribution}{Stationary Markov chain with Bern(p) marginal distribution}}
The literature has extensively discussed the statistical inference for finite Markov chains. When transition probabilities of the Markov chain are expressed as $p_{ij}(\boldsymbol{\theta})$, with ${\boldsymbol{\theta}}=(\theta_{{1}},\cdots,\theta{r})$, the log-likelihood is $\sum_{ij}n_{ij}\log p_{ij}(\boldsymbol{\theta})$ (see Billingsley (1961a) \cite{billingsley1961statistical} and Billingsley (1961b) \cite{billingsley1961methods}). If the $r$ components of $\boldsymbol{\theta}$ are all real, then the maximum likelihood equations are \begin{align}\label{MLE equations}\sum_{ij}\frac{n_{ij}}{p_{ij}(\boldsymbol{\theta})}\frac{\partial p_{ij}(\boldsymbol{\theta})}{\partial \theta_u}=0,\ u=0,1,\cdots,r.   \end{align} 
Considering the transition matrix (\ref{tmatrples}) for instance, the likelihood function of the sample is
	\begin{align}\label{likelihood 1}
		L(\boldsymbol{\theta})=p^{x_{0}}(1-p)^{1-x_{0}}  \left(\frac{ap+1-2p}{1-p}\right)^{n_{_{00}}}\left(\frac{p(1-a)}{1-p}\right)^{n_{_{01}}}(1-a)^{n_{_{10}}}a^{n_{_{11}}}
	\end{align} and the log-likelihood is given by\begin{equation}\label{loglikelihood}
		\begin{split}
			\mathcal{L}(\boldsymbol{\theta})=\log L(\boldsymbol{\theta})=&x_{0}\log(p)+(1-x_{0})\log(1-p)+n_{00}\log(ap+1-2p)-\\& n_{00}\log(1-p)+n_{01}\log(p)+ n_{01}\log(1-a)-n_{01}\log(1-p)\\& \qquad + n_{10}\log(1-a)+n_{11}\log(a).
		\end{split}
	\end{equation}
Differentiating (\ref{loglikelihood}) with respect to $a$ and $p$ gives 
	\begin{equation}\label{likeqns}
		\begin{cases}
			&\frac{\partial \mathcal{L}(\boldsymbol{\theta})}{\partial a}=\frac{n_{_{00}}p}{ap+1-2p}-\frac{n_{_{01}}}{1-a}-\frac{n_{_{10}}}{1-a}+\frac{n_{11}}{a}=0\\
			&\frac{\partial \mathcal{L}(\boldsymbol{\theta})}{\partial p}=\frac{x_{0}}{p}-\frac{1-x_{0}}{1-p}+\frac{n_{_{00}}(a-2)}{ap+1-2p}+\frac{n_{_{00}}}{1-p}+\frac{n_{_{01}}}{p}+\frac{n_{_{01}}}{1-p}=0.
		\end{cases}
	\end{equation}
	The equations in (\ref{likeqns}) can be rewritten in the form\begin{align}\label{likeqns simplified}
		\begin{cases}
			&npa^2-[(2n-n_{_{00}}+n_{_{11}})p-n+n_{_{00}}]a+n_{_{11}}(2p-1)=0\\
			& a[x_0p-p^2+n_{_{00}}p+n_{_{01}}p]+x_0-2x_0p-p+2p^2-n_{_{00}}p+n_{_{01}}-2n_{_{01}}p =0.
		\end{cases}
	\end{align}
	Solving for $a$ in the second equation of (\ref{likeqns simplified}) gives \begin{align}\label{maxlik solved for a}
		\hat{a}=\frac{p(\lambda_1-2p)-\lambda_2}{p(\lambda_3-p)},  \end{align} where $\lambda_1=2x_0+1+n_{_{00}}+2n_{_{01}},\lambda_2=x_0+n_{_{01}}$ and $\lambda_3=x_0+n_{_{00}}+n_{_{01}}.$
	Substituting (\ref{maxlik solved for a}) for $a$ in the first equation of (\ref{likeqns simplified}) and using $\lambda_4=2n-n_{_{00}}+n_{_{11}},$ and $\lambda_5=n-n_{_{00}}$  we obtain \begin{align}\label{values of a and p}    (4n-2\lambda_4&+2n_{_{11}})p^4+(2\lambda_3\lambda_4-4n\lambda_1+\lambda_1\lambda_4+2\lambda_5-n_{_{11}}-4n_{_{11}}\lambda_3)p^3+\nonumber\\&(n\lambda_1^2+4n\lambda_2-\lambda_1\lambda_3\lambda_4-\lambda_2\lambda_4-2\lambda_3\lambda_5-\lambda_1\lambda_5+2n_{_{11}}\lambda_3+2n_{_{11}}\lambda_3^2)p^2+\nonumber\\&(\lambda_1\lambda_3\lambda_5-2n\lambda_1\lambda_2+\lambda_2\lambda_3\lambda_4+\lambda_2\lambda_5-n_{_{11}}\lambda_3^2)p+(n\lambda_2^2-\lambda_2\lambda_3\lambda_5)=0,\end{align}
	There is no closed-form solution for equation (\ref{values of a and p}) but numerical solutions can be obtained.
	
	The sufficient conditions for MLEs and their asymptotic properties are specified in Condition 5.1 of Billingsley (1961a) \cite{billingsley1961statistical}, as presented below.
	
	\begin{condition}\label{condi1}
		\textit{Let $X=(X_t)_{t=0,\cdots,n}$ be a Markov chain with finite state space $S=\lbrace 0,\cdots,s\rbrace$. Suppose $P=(p_{ij}(\boldsymbol{\theta}))$ is the $s\times s$ stochastic matrix, where ${\boldsymbol{\theta}}=(\theta_{_{1}},\cdots,\theta_{r})$ ranges over an open subset $\Theta$ of $\mathbb{R}^r$. Suppose that the following conditions are satisfied:
			\begin{itemize}
				\item[i.]   The set  $E$ of $(i,j)$ such that $p_{ij}(\boldsymbol{\theta})>0$ is independent of $\boldsymbol{\theta}$ and each $p_{ij}(\theta)$ has continuous partial derivatives of third order throughout $\Theta.$
				\item[ii.] The $d\times r$ matrix $D=(\frac{\partial p_{ij}(\boldsymbol{\theta})}{\partial \theta_{k}})$, $k=1,\cdots, r$, where $d$ is the number of elements in $E$, has rank $r$ throughout $\Theta$.
				\item[iii.] For each $\boldsymbol{\theta}\in\Theta$, there is only one ergodic set and no transient states.
			\end{itemize}
			Then there exists $\boldsymbol{\hat{\theta}}\in \Theta$ such that, with probability tending to $1$, $\boldsymbol{\hat{\theta}}$ is a solution of system (\ref{MLE equations}) and converges in probability to the true value $\boldsymbol{\theta}$. Additionally, $\sqrt{(n+1)}(\boldsymbol{\hat{\theta}}-\boldsymbol{\theta})\xrightarrow{d}N_{r}(\boldsymbol{0},\boldsymbol{\Sigma^{-1}})$, where $\boldsymbol{\Sigma^{-1}}$ is the inverse of matrix $\boldsymbol{\Sigma}$ with entries defined by
			\begin{equation}\label{vars}
				\sigma_{uv}(\boldsymbol{\theta})=E_{\boldsymbol{\theta}}\left[\frac{\partial \log p_{x_{1}x_{2}}(\boldsymbol{\theta})}{\partial \theta_{u}}\frac{\partial \log p_{x_{1}x_{2}}(\boldsymbol{\theta})}{\partial \theta_{v}} \right],\qquad u,v=1,\cdots, r.
			\end{equation}.}
	\end{condition}
	This condition can be used to verify the following theorem:
	\begin{theorem}\label{Estimators for discrete}
		Let $X=(X_t)_{t=0,\cdots,n}$ be a stationary Markov chain generated by copula (\ref{our model}) and Bernoulli($p$) marginal distribution. Then there exists a consistent MLE, $\boldsymbol{\hat{\theta}}=(\hat{a},\hat{p})\ \text{of }\ \boldsymbol{\theta}=(a,p)$ such that   $ \sqrt{(n+1)}(\hat{\boldsymbol{\theta}}-\boldsymbol{\theta})\xrightarrow{d}N_{2}(\boldsymbol{0},\boldsymbol{\Sigma^{-1}}),$ where  \begin{align}\label{sigmainverse}
			\Sigma^{-1}=\begin{cases}
				\begin{pmatrix}
					\frac{a(1-a)}{p}& a(1-p)\\
					a(1-p) & \frac{p(1-p)(a+1-2p)}{1-a}
				\end{pmatrix}, \ \text{when}\ p<1/2\\  \begin{pmatrix}
					\frac{a(1-a)}{1-p}& -ap\\
					-ap & \frac{p(1-p)(2p-1+a)}{1-a}
				\end{pmatrix},\ \text{when}\ p> 1/2.  
			\end{cases}
		\end{align} 
	\end{theorem}
	Note that $p=1/2$ does not fall into either of these two cases because Condition (\ref{condi1}) applies only to interior points of the range of $(a,p)$.
	
	\begin{proof}[Proof of Theorem (\ref{Estimators for discrete})]
		The proof of Theorem (\ref{Estimators for discrete}) follows by first verifying that the requirements outlined in Condition (\ref{condi1}) are met for each of the cases, $p<1/2$ and $p> 1/2$.   The entries of the information matrices, $\boldsymbol{\Sigma}$, are obtained by computing the expectations given in Equation (\ref{vars}). Inverting $\boldsymbol{\Sigma}$, gives (\ref{sigmainverse}). The complete proof can be found in Muia (2024) \cite{muia2024dependence}.
	\end{proof}
	
	When $p=1/2$, the stationary Markov chain based on copula (\ref{our model}) and \\ Bernoulli($1/2$) margins has a transition matrix of the form \begin{equation}\label{nonstationarymatrix}
		P=  \begin{pmatrix}
			a & 1-a\\
			1-a& a
		\end{pmatrix}.
	\end{equation} 
	In this case we are left with $a$ to estimate. The likelihood function for the chain is given by 
	\begin{align*}
		L(\textbf{x},a)=\frac{1}{2}a^{n_{00}+n_{11}}(1-a)^{n-(n_{00}+n_{11})},
	\end{align*} 
	and the log-likelihood is expressed as
	\begin{align}\label{loglike3}
		\mathcal{L}(\textbf{x};a)=(n_{00}+n_{11})log(a)+(n-(n_{00}+n_{11}))log(1-a).   
	\end{align}
	The following proposition provides the MLE of $a$ and its asymptotic properties in this case.
	\begin{prop}\label{prop p=1/2}
		Let $(X_t)_{t=0,\cdots, n}$ be a stationary Markov chain generated by copula (\ref{our model}) and Bernoulli($p$) marginal distribution. When $p=1/2$, the parameter $a$ has  MLE $\hat{a}=\frac{n_{00}+n_{11}}{n}$. Moreover, $\frac{\sqrt{(n+1)}(\hat{a}-a)}{\sqrt{a(1-a)}}\xrightarrow{d}N(0,1).$
	\end{prop}
	\begin{proof}
		The proof of Proposition (\ref{prop p=1/2}) follows from verifying the requirements of Condition (\ref{condi1}). The variance of $a$ is computed using formula (\ref{vars}) with $r=1$.    
	\end{proof}

	Klotz showed that the MLE, $\boldsymbol{\hat{\theta}}=(\hat{a},\hat{p})$, is asymptotically equivalent to  $(a^{*}(\Bar{p}),\Bar{p})$, where $a^{*}(\Bar{p})$ is the MLE of $a$ when $p$ is replaced by $\Bar{p}=1/(n+1)\sum_{i=0}^n X_i.$ The proof relies on the fact that $\hat{p}$ and $\bar{p}$ have the same asymptotic variance.  Lindqvist (1978) \cite{lindqvist1978note} considers a general case of non-Markovian dependence and shows that the MLE, $\hat{p}$ is asymptotically equivalent to $\Bar{p}$ given that the process is ergodic.
	
	In this work, different asymptotically equivalent estimators for $p$ are discussed. Since the Markov chain under consideration is stationary, it follows that $P(X_t)=p\ \forall\ t$. The  sample proportion, $\Bar{p}=1/(n+1)\sum_{i=0}^n X_i$, is an unbiased estimator for $p$ (see Bedrick and Aragon (1989) \cite{bedrick1989approximate}, Klotz (1973) \cite{klotz1973statistical}) or Lindqvist (1978) \cite{lindqvist1978note}). 
	
	\begin{prop}
		Let $X=(X_t)_{t=0,\cdots,n}$ be a stationary Markov chain generated by  copula (\ref{our model}) and Bernoulli($p$) marginal distribution.
		The asymptotic variance of the sample proportion, $\bar{p},$ is the same as that of $\hat{p}$.
	\end{prop}
	\begin{proof}
		For the Markov chain generated when $p<1/2$
		{\small \begin{align*}
				var(\Bar{p})&=var\left(\frac{1}{n+1}\sum_{t=0}^nX_t\right)\\
				&=\frac{1}{(n+1)^2}\left[(n+1)var(X_0)+2\sum_{t=0}^{n-1}\sum_{k=1}^{n-t}cov(X_t,X_{t+k}) \right]\\
				&=\frac{1}{(n+1)^2}\left[(n+1)p(1-p) +2p(1-p)\sum_{t=0}^{n-1}\sum_{k=1}^{n-t}\left( \frac{a-p}{1-p}\right)^k\right]
			\end{align*}
			Applying the properties of a geometric sum twice and simplifying, we obtain
			\begin{align*}
				var(\Bar{p})=\frac{p(1-p)}{n+1}\left[1 +2\left(\frac{a-p}{1-a}\right) \left(\frac{n}{n+1}+\frac{1}{n+1} \left( \frac{a-p}{1-a}\right)\left(\left( \frac{a-p}{1-p}\right)^n-1\right)\right)\right].
		\end{align*} }
		Therefore, for large $n$,
		{\small \begin{align}\label{varbarX1}
				var(\Bar{p})\approx \frac{p(1-p)}{n+1}\left[ 1 +\frac{2a-2p}{1-a}\right]+O\left(\frac{1}{n} \right)=\frac{p(1-p)}{n+1}\left[\frac{a+1-2p}{1-a}\right]+O\left(\frac{1}{n} \right)\end{align} }
		Similarly, for $p\geq 1/2$, 
		{\small  \begin{align}\label{varbarX2}
				var(\Bar{p})\approx \frac{p(1-p)}{n+1}\left[ 1 +\frac{2a+2p-2}{1-a}\right]+O\left(\frac{1}{n} \right)=\frac{p(1-p)}{n+1}\left[\frac{2p-1+a}{1-a}\right]+O\left(\frac{1}{n} \right)\end{align} }
		From (\ref{varbarX1}) and (\ref{varbarX2}) we note that the asymptotic variance of the sample mean is the same as that of $\hat{p}$ (refer to Theorem (\ref{Estimators for discrete}) for the asymptotic variance of $\hat{p}$).
	\end{proof}

	The following result provides the asymptotic properties of the sample mean for the stationary Markov chain based on copula (\ref{our model}).
	
	\begin{theorem}\label{clt of xbar}
		Let $X=(X_t)_{t=0,\cdots,n}$ be a stationary Markov chain generated by  copula (\ref{our model}) and Bernoulli($p$) marginal distribution. Let $ \Bar{p}$ be the sample mean. It follows that 
		\begin{align*}
			\sqrt{(n+1)}(\Bar{p}-p)\xrightarrow{d} N(0,\sigma^2),
		\end{align*}
		where $\sigma^2=p(1-p)(1+a-2p)/(1-a)$, $p<1/2$  and $\sigma^2=p(1-p)(a+2p-1)/(1-a)$, $p\geq 1/2$. 
	\end{theorem}
	
	\begin{proof}
		The proof of Theorem (\ref{clt of xbar})  follows from Theorem (18.5.2) of Ibragimov and Linnik (1975) \cite{ibragimov1975independent} and the fact that the Markov chain is uniformly mixing.  Using steps similar to those used in the proof of Theorem (\ref{discrete}), the $\phi$-mixing coefficient is found to be \begin{align*}\phi(n)=\begin{cases}
				(1-p)\left|\frac{a-p}{1-p}\right|^n, &p<1/2\\
				p\left|\frac{a+p-1}{p}\right|^n, &p\geq1/2.
		\end{cases}\end{align*} Let $p<1/2$. To check Condition (18.5.8) of Theorem (18.5.2) of Ibragimov and Linnik (1975) \cite{ibragimov1975independent}, we see that \begin{align}\label{geometricseries1}
			\sum_{n=0}^\infty \lbrace \phi(n)\rbrace^{1/2}=\sum_{n=0}^\infty\sqrt{(1-p)}\left(\left|\frac{a-p}{1-p}\right|^{1/2}\right)^n,
		\end{align}
		converges since $\left|(a-p)/(1-p)\right|^{1/2}<1$. The proof for $p\geq 1/2$ is identical. 
		The variance of $\bar{X_n}$ is given by the formulae (\ref{varbarX1}) and (\ref{varbarX2}). Hence Condition (18.5.9) of Theorem (18.5.2) of Ibragimov and Linnik (1975) \cite{ibragimov1975independent}  satisfied.
		
	\end{proof}

\section{Hypothesis testing}\label{test for independence}
\subsection{Test of Independence}
Different tests have been used in the literature to assess the independence of observations in a sequence. These include the $\chi^2$ goodness-of-fit tests and the Likelihood Ratio Test (LRT), as discussed by Anderson and Goodman (1957) \cite{anderson1957statistical}, Goodman (1959) \cite{goodman1959some}, Billingsley (1961a) \cite{billingsley1961statistical}, and others. 

Having $p_{ij} = p_j$ for all $j$ is sufficient for independence in a finite state-space Markov chain. The following proposition outlines conditions for independence in the Markov chain based on copula (\ref{our model}) and Bernoulli($p$) marginal distributions.

\begin{prop}\label{theorem for independence}
	Let $(X_t)_{t=0,\cdots,n}$ be a stationary Markov chain with copula (\ref{our model}) and Bernoulli($p$) marginal distribution. Then $X_0,\cdots,X_n$ are independent when $a=p<1/2$ or $a=1-p\leq 1/2$.
\end{prop}
\begin{proof}
	For independent observations in a sequence, \(p_{ij} = p_j\) for all states \(i,j\), where \(p_j\) is the probability of outcome \(j\). When \(p<1/2\), conditions for independence are given by \((ap+1-2p)/(1-p)= 1-a\) and \(p(1-a)/(1-p)=a\), leading to \(a=p\) in both cases. For \(p\geq 1/2\), independence holds if \(a=(1-p)(1-a)/p\) and \(1-a=(a(1-p)+2p-1)/p\), simplifying to \(a=1-p\) in each case.
\end{proof}

The likelihood ratio approach for testing the hypothesis of independence is obtained as follows: For $p<1/2$, the likelihood function is expressed as shown in (\ref{likelihood 1}). The hypotheses under consideration for the test are denoted as $H_0:\theta\in\omega$ versus $H_1:\theta\in\Theta\cap\omega^c$, where
$\Theta = \lbrace\theta= (a,p), 0<a<1, 0<p<1/2\rbrace$ and
$\omega = \lbrace\theta= (p,p), 0<p<1/2\rbrace\subset\Theta.$
Let $\hat{\theta}$ be the MLE under $H_a$ and $\hat{\theta}_0$ the MLE under $H_0$. 
The LR test statistic is then given by \begin{align}
	\Lambda=\frac{\max\limits_{\theta\in\omega}L(\theta)}{\max\limits_{\theta\in\Theta}L(\theta)}=\frac{L(\hat{\theta}_0)}{L(\hat{\theta})}.
\end{align}
The likelihood under $H_0$ is expressed as:\begin{align}\label{simpler model lik}
	L(\omega)=p^{x_0+n_{_{01}}+n_{_{11}}}(1-p)^{1-x_0+n_{_{00}}+n_{_{10}}}.
\end{align}
The  function (\ref{simpler model lik}) represents the likelihood of an i.i.d Bern($p$) sequence. Note that $x_0+n_{_{01}}+n_{_{11}}=X=\text{number of ones}.$  The MLE under $H_0$ is\begin{align*}
	\hat{p}=\frac{x_0+n_{_{01}}+n_{_{11}}}{n+1}=\frac{X}{n+1}\ \text{and}\ \hat{\theta}_0=(\hat{p},\hat{p}).
\end{align*}
The likelihood function evaluated at the MLE under $H_0$ simplifies to \begin{align*}
	L(\hat{\omega})=\left(\frac{X}{n+1}\right)^{X}\left(1-\frac{X}{n+1}\right)^{n+1-X}.
\end{align*}
The LRT statistic is therefore \begin{align}\label{LRT ples}
	\Lambda_1   = \frac{\left(\frac{X}{n+1}\right)^{X}\left(1-\frac{X}{n+1}\right)^{n+1-X}}{\hat{p}^{x_{_0}}(1 - \hat{p})^{1-x_{_0}}\left( \frac{\hat{a}\hat{p}+1-2\hat{p}}{1-\hat{p}} \right)^{n_{_{00}}}\left( \frac{\hat{p}(1-\hat{a})}{1-\hat{p}} \right)^{n_{_{01}}}(1-\hat{a})^{n_{_{10}}} \hat{a}^{n_{_{11}}}}.
\end{align}
For a given significance level $\alpha$, reject $H_0$ in favor of $H_1$ if $\Lambda_1<c$, where $c$ is such that $\alpha=P_{H_0}[\Lambda_1\leq c]$. 

A result due to Wilks (1938) \cite{wilks1938large} shows that under suitable regularity conditions, if $\Theta$ is a $r$-dimensional space and $\omega$ is a $k$-dimensional space, then $-2\log\Lambda_1$ has an asymptotically $\chi^2$ distribution with $r-k$ degrees of freedom. The fulfillment of these desired conditions has been confirmed for our model by demonstrating that the transition probabilities satisfy the requirements outlined in Condition (\ref{condi1}). The asymptotic distribution of $-2\log\Lambda_1$ is $\chi^2$ with 1 degree of freedom. The decision rule is to reject $H_0$ in favor of $H_1$ if $-2\log\Lambda_1\geq\chi^2_{\alpha}(1)$. 

In a similar manner, for $p> 1/2$, define $    \Theta=\lbrace \theta=(a,p): 0<a<1,1/2< p<1\rbrace$ and $\omega=\lbrace \theta=(1-p,p): 1/2< p<1\rbrace.$ The distribution doesn't change under $H_0$, but changes under $H_a$. 
The LRT statistic becomes\begin{align}\label{LRT pgeq}
	\Lambda_2   = \frac{\left(\frac{X}{n+1}\right)^{X}\left(1-\frac{X}{n+1}\right)^{n+1-X}}{\hat{p}^{x_{_0}}(1 - \hat{p})^{1-x_{_0}}\left(\frac{(1-\hat{p})(1-\hat{a})}{\hat{p}}\right)^{n_{_{10}}} \left( \frac{\hat{a}(1-\hat{p})+2\hat{p}-1}{\hat{p}}\right)^{n_{_{11}}}\left(1-\hat{a} \right)^{n_{_{01}}}\hat{a} ^{n_{_{00}}}}.
\end{align}
The limiting distribution of $-2\log\Lambda_2$ approaches a chi-squared distribution with one degree of freedom as $n\longrightarrow\infty$. The level $\alpha$ decision rule is to reject\ $H_0$\ in favor of $H_1$ if  $-2\log\Lambda_2\geq\chi^2_{\alpha}(1)$.

\section{Simulation Study}\label{Simulations Study}
\subsection{MLEs of a and p}
Tables \ref{MLE1} and \ref{MLE2} below present the MLEs of $a$ and $p$ for various sample sizes and selected initial value pairs. Few pairs of $a$ and $p$ were selected for this study to avoid very lengthy tables. We used different combinations of true values for $a$ and $p$ and constructed 95\% confidence intervals accordingly. Utilizing the asymptotic normality described in Theorem (\ref{Estimators for discrete}), we derived confidence intervals for $a$ and $p$ based on the standard normal distribution.Furthermore, we generated 400 replications for each sequence of length $n$, which were used to calculate the coverage probabilities (CPs) of the true values by the confidence intervals, as well as the mean lengths of the confidence intervals (CIML). Note that we did not include the numerical values of the estimators because $400$ numerical values were computed for each estimator in each case, and there is no criterion to choose only one representative from the $400$. Relatively narrow confidence intervals reflect a higher level of precision in estimation, while wider confidence intervals indicate greater uncertainty in estimation.

In Table \ref{MLE1}, the confidence intervals for $a$ decrease as $p$ increases, reflecting reduced uncertainty in estimating $a$ with higher values of $p<0.5$. This reduction is attributed to the decreasing variance of $\hat{a}$ as $\hat{p}<1/2$ increases. The accuracy in estimation improves with sample size. Moreover, as the sample size grows, the coverage probability for $a$ becomes more precise, indicating improved accuracy in capturing the true value of $a$ by the confidence intervals.

\begin{table}
	\centering
	\caption{The average length of $400,\ 95\%$ confidence intervals and the coverage probabilities for the true values of $a$ and $p$ when $p<1/2$.}
	\label{MLE1}
	{\scriptsize \begin{tabular}{@{}c|ccccc@{}}
			\toprule
			$n$&Initial parameters & CIML for `a' & CIML for `p' & CP for `a' & CP for `p' \\ \midrule \midrule
			\multirow{11}{*}{499}&$a=0.1, p=0.1$ & 0.1586 & 0.0523 & 0.8925 & 0.9425 \\ 
			&$a=0.1, p=0.3$ & 0.0950 & 0.0598 & 0.9300 & 0.9500 \\
			&$a=0.1, p=0.4$ & 0.0821 & 0.0494 & 0.9475 & 0.9500 \\
			&$a=0.2, p=0.1$ & 0.2172 & 0.0584 & 0.9050 & 0.9425 \\
			&$a=0.2, p=0.3$ & 0.1273 & 0.0695 & 0.9475 & 0.9575 \\
			&$a=0.2, p=0.4$ & 0.1100 & 0.0605 & 0.9475 & 0.9525 \\
			&$a=0.7, p=0.1$ & 0.2671 & 0.1153 & 0.9675 & 0.9050 \\
			&$a=0.7, p=0.3$ & 0.1479 & 0.1530 & 0.9650 & 0.9625 \\
			&$a=0.7, p=0.4$ & 0.1275 & 0.1490 & 0.9500 & 0.9475 \\
			&$a=0.9, p=0.1$ & 0.2237 & 0.2000 & 0.9500 & 0.8350 \\
			&$a=0.9, p=0.3$ & 0.1014 & 0.2830 & 0.9675 & 0.9175 \\
			&$a=0.9, p=0.4$ & 0.0869 & 0.2797 & 0.9650 & 0.8825 \\ \midrule
			\multirow{11}{*}{999}&$a=0.1, p=0.1$ & 0.1146 & 0.0370 & 0.9125 & 0.9325 \\ 
			&$a=0.1, p=0.3$ & 0.0674 & 0.0423 & 0.9525 & 0.9500 \\
			&$a=0.1, p=0.4$ & 0.0585 & 0.0351 & 0.9275 & 0.9350 \\
			&$a=0.2, p=0.1$ & 0.1554 & 0.0414 & 0.9300 & 0.9200 \\
			&$a=0.2, p=0.3$ & 0.0902 & 0.0491 & 0.9475 & 0.9350 \\
			&$a=0.2, p=0.4$ & 0.0782 & 0.0429 & 0.9500 & 0.9400 \\
			&$a=0.7, p=0.1$ & 0.1854 & 0.0818 & 0.9600 & 0.9225 \\
			&$a=0.7, p=0.3$ & 0.1042 & 0.1084 & 0.9550 & 0.9550 \\
			&$a=0.7, p=0.4$ & 0.0899 & 0.1050 & 0.9350 & 0.9575 \\
			&$a=0.9, p=0.1$ & 0.1373 & 0.1446 & 0.9625 & 0.8825 \\
			&$a=0.9, p=0.3$ & 0.0693 & 0.2035 & 0.9500 & 0.9275 \\
			&$a=0.9, p=0.4$ & 0.0597 & 0.1997 & 0.9550 & 0.9600 \\ \midrule
			\multirow{11}{*}{4999}&$a=0.1, p=0.1$ & 0.0524 & 0.0166 & 0.9450 & 0.9600 \\
			&$a=0.1, p=0.3$ & 0.0303 & 0.0189 & 0.9575 & 0.9525 \\
			&$a=0.1, p=0.4$ & 0.0263 & 0.0157 & 0.9475 & 0.9500 \\
			&$a=0.2, p=0.1$ & 0.0700 & 0.0186 & 0.9600 & 0.9525 \\
			&$a=0.2, p=0.3$ & 0.0404 & 0.0220 & 0.9500 & 0.9475 \\
			&$a=0.2, p=0.4$ & 0.0350 & 0.0192 & 0.9600 & 0.9325 \\
			&$a=0.7, p=0.1$ & 0.0809 & 0.0370 & 0.9600 & 0.9325 \\
			&$a=0.7, p=0.3$ & 0.0465 & 0.0486 & 0.9600 & 0.9650 \\
			&$a=0.7, p=0.4$ & 0.0402 & 0.0471 & 0.9675 & 0.9775 \\
			&$a=0.9, p=0.1$ & 0.0541 & 0.0678 & 0.9300 & 0.9425 \\
			&$a=0.9, p=0.3$ & 0.0306 & 0.0913 & 0.9425 & 0.9525 \\
			&$a=0.9, p=0.4$ & 0.0264 & 0.0899 & 0.9550 & 0.9525 \\ \bottomrule
	\end{tabular} }
\end{table}

Table \ref{MLE2} shows that confidence intervals lengthen as $p$ increases due to the increasing variance of $\hat{a}$ for $\hat{p}>1/2$. The coverage probabilities for both $a$ and $p$ vary with different combinations of initial values and sample sizes, underscoring the significance of careful parameter selection in statistical inference. Nevertheless, as sample size increases, both coverage probabilities tend to converge towards the desired nominal value of $95\%$.

Furthermore, the results presented in Tables \ref{MLE1} and \ref{MLE2} demonstrate a noteworthy symmetry in the CIMLs when comparing parameter pairs equidistant from \( p = 0.5 \) for fixed \( a \). That is, for each fixed value of \( a \), and for \( p \) values equidistant from \( p = 0.5 \), the CIMLs for \( a \) and \( p \) are the same. Figure \ref{symmetry} below shows this symmetry for $n=9999$. 

The observed behavior of the mean length confidence intervals (CIML) for \(a\) and \(p\) can be explained by considering the variance-covariance matrix \(\Sigma^{-1}\) from Theorem \ref{Estimators for discrete}, which governs the precision of the maximum likelihood estimators (MLEs) \(\hat{a}\) and \(\hat{p}\). When \(p\) is closer to 0.5, the off-diagonal elements (representing the covariance between \(a\) and \(p\)) and the denominators in the variance terms become more influential. As \(p\) increases from 0.1 to 0.4 (for \(p < 0.5\)), the CIML for \(a\) decreases, indicating reduced uncertainty in estimating \(a\) due to lower variance. For \(p > 0.5\), as \(p\) moves further away from 0.5, the CIML for \(a\) increases again because of increased variance terms and covariance between \(a\) and \(p\). The symmetry in the CIML for \(a\) around \(p = 0.5\) (e.g. similar CIML for \(p = 0.1\) and \(p = 0.9\)) suggests that confidence intervals widen as \(p\) deviates from 0.5 in either direction, reflecting greater difficulty in estimating \(a\) at extreme values of \(p\).

\begin{table}
	\centering
	\caption{The average length of $400,\ 95\%$ confidence intervals and the coverage probabilities for the true values of $a$ and $p$ when $p> 1/2$.}
	\label{MLE2}
	{\scriptsize  \begin{tabular}{c|cccccc}
			\toprule
			$n$ & $a$ & $p$ &  CIML for `a' & CIML for `p' & CP for `a' & CP for `p' \\ \midrule\midrule
			\multirow{12}{*}{499}&0.1 & 0.6  & 0.0821       & 0.0494       & 0.9475   & 0.9500 \\
			&0.1 & 0.7 & 0.0950       & 0.0598       & 0.9300   & 0.9500 \\
			&0.1 & 0.9 &  0.1593       & 0.0524       & 0.8925   & 0.9425 \\
			&0.3 & 0.6 &  0.1267       & 0.0726       & 0.9500    & 0.9575\\
			&0.3 & 0.7 & 0.1460       & 0.0803       & 0.9275   & 0.9550 \\
			&0.3 & 0.9 &  0.2514       & 0.0658       & 0.9325    & 0.9550\\
			&0.7 & 0.6 &  0.1275       & 0.1484       & 0.9475   &  0.9475\\
			&0.7 & 0.7 &  0.1479       & 0.1530       & 0.9650   &0.9625  \\
			&0.7 & 0.9 &  0.2671       & 0.1153       & 0.9675   & 0.905 \\
			&0.9 & 0.6 &  0.0868       & 0.2797       & 0.9600   &  0.9625\\
			&0.9 & 0.7 &  0.1014       & 0.2830       & 0.9675   &  0.9175\\
			&0.9 & 0.9 &  0.2237       & 0.2000       & 0.9500   &  0.8350\\ \midrule
			\multirow{12}{*}{999}&0.1 & 0.6  & 0.0585       & 0.0351       & 0.9275   &  0.9350\\
			&0.1 & 0.7 &  0.0674       & 0.0423       & 0.9525    & 0.9500\\
			&0.1 & 0.9 & 0.1146       & 0.0370       & 0.9125    & 0.9325\\
			&0.3 & 0.6 &  0.0896       & 0.0512       & 0.9500   &  0.9450\\
			&0.3 & 0.7 & 0.1034       & 0.0567       & 0.9450   &  0.9325\\
			&0.3 & 0.9 &  0.1787       & 0.0464       & 0.9475    & 0.9350\\
			&0.7 & 0.6 &  0.0899       & 0.1050       & 0.9350    & 0.9575\\
			&0.7 & 0.7 &  0.1042       & 0.1084       & 0.9550   &  0.9550\\
			&0.7 & 0.9 & 0.1854       & 0.0818       & 0.9600   &  0.9225\\
			&0.9 & 0.6 &  0.0596       & 0.1997       & 0.9550   &  0.9550\\
			&0.9 & 0.7 &  0.0693       & 0.2035       & 0.9500   &  0.9275\\
			&0.9 & 0.9 & 0.1373       & 0.1446       & 0.9625   &  0.8825\\ \midrule
			\multirow{12}{*}{4999}&0.1 & 0.6 &  0.0263     & 0.0157       & 0.9475   &  0.9500\\
			&0.1 & 0.7 &  0.0303     & 0.0189       & 0.9575   &  0.9525\\
			&0.1 & 0.9 &  0.0524     & 0.0166       & 0.9450   &  0.9600\\
			&0.3 & 0.6 &  0.0401     & 0.0229       & 0.9300   &  0.9450\\
			&0.3 & 0.7 &  0.0463     & 0.0254       & 0.9350   &  0.9425\\
			&0.3 & 0.9 &  0.0803     & 0.0208       & 0.9525   &  0.9325\\
			&0.7 & 0.6 &  0.0402     & 0.0471       & 0.9675   &  0.9775\\
			&0.7 & 0.7 &  0.0465     & 0.0486       & 0.9600    & 0.9650\\
			&0.7 & 0.9 &  0.0809     & 0.0370       & 0.9600   &  0.9325\\
			&0.9 & 0.6 &  0.0264     & 0.0899       & 0.9550   &  0.9525\\
			&0.9 & 0.7 &  0.0306     & 0.0913       & 0.9425   &  0.9525\\
			&0.9 & 0.9 &  0.0541     & 0.0678       & 0.9300   &  0.9325\\ \bottomrule
	\end{tabular} }
\end{table} 

For \(p\), the mean length of the confidence intervals is relatively stable across different values, with small fluctuations likely due to the interaction between \(a\) and \(p\) in the variance terms of \(\Sigma^{-1}\). For small \(a\) (e.g., \(a = 0.1\)), the CIML for \(p\) is slightly lower when \(p\) is near 0.1 and 0.9 compared to intermediate values, indicating higher certainty in extreme probability scenarios. Given the large sample size (\(n = 9999\)), the confidence intervals are generally narrow, as the variance of the MLEs decreases with increasing sample size, leading to more precise estimates. 

Overall, the behavior of the CIML for \(a\) and \(p\) reflects the interplay between the values of \(a\) and \(p\), their interaction in the variance-covariance structure, and the large sample size ensuring narrow confidence intervals. The symmetry around \(p = 0.5\) and the increasing interval lengths as \(p\) approaches the extremes illustrate the underlying statistical properties of the MLEs and their variances.

This symmetry reflects the underlying properties of the mixture copula model. The behavior suggests that the model treats the transition probabilities in a symmetric manner around \( p = 0.5 \), ensuring that the statistical characteristics such as interval lengths and coverage probabilities remain stable for \( p \) values that are symmetrically opposite about 0.5. This property enhances the robustness of the parameter estimation process within the specified range of \( p \).  

	\begin{figure}[ht]
		\centering
		%\begin{mdframed}
		\fbox{\includegraphics[scale=0.5]{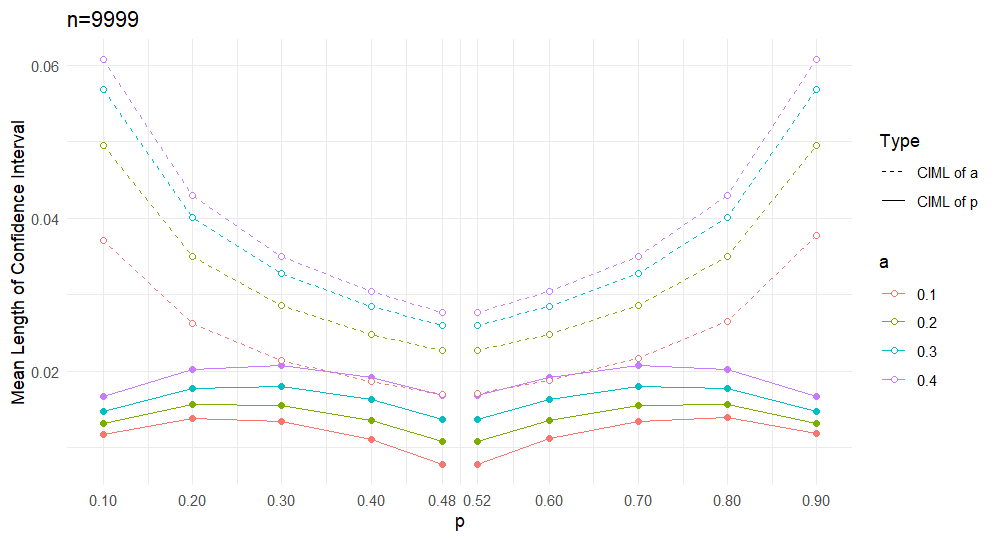}}
		%\end{mdframed}
		\caption{Symmetry in Mean Lengths of Confidence Intervals for `p' about $p=0.5$}
		\label{symmetry}
	\end{figure}

\subsection{Test for independence}
\subsubsection{Likelihood ratio test for independence}
Tables \ref{chi_squared_results1} and \ref{chi_squared_results2} below present the results of the likelihood ratio test for independence in the Markov chain generated by copula model (\ref{our model}) and Bernoulli($p$) marginal distribution. The tables provide separate results for $p<1/2$ and $p>1/2$, respectively.

The Markov chain is generated using true values of $a$ and $p$, and subsequently, a chi-squared independence test is performed based on the test statistics $-2\log\Lambda_1$ and $-2\log\Lambda_2$, where $\Lambda_1$ and $\Lambda_2$ are given by formulae (\ref{LRT ples}) and (\ref{LRT pgeq}) respectively. A significance level of $\alpha = 0.05$ is chosen for the test, implying that the null hypothesis of independence will be rejected if  $-2\log\Lambda_i\geq 3.842, i=1,2$. Conversely, if $-2\log\Lambda_i< 3.842$, the null hypothesis of independence is not rejected.

The chain is simulated for a total of $n + 1$ steps, where $n$ is set to $9999$. The iterative process for Table \ref{chi_squared_results1} has true values of $a$ ranging from $0.1$ to $0.9$ with a step of $0.1$ and $p < 1/2$ ranging from $0.1$ to $0.4$ with a step of $0.1$. On the other hand, Table \ref{chi_squared_results2} has true values of $a$ ranging from $0.1$ to $0.9$ with a step of $0.1$ and $p >1/2$ ranging from $0.6$ to $0.9$ with a step of $0.1$.

\begin{table}[htbp]
	\centering
	%\caption{Likelihood Ratio Test Results and KS Distance values for different combinations of $a$ and $p$}
	\begin{minipage}{0.5\textwidth}
		\centering
		\setlength{\tabcolsep}{6pt} % Adjust column separation to match both tables
		\renewcommand{\arraystretch}{1.2} % Adjust row height
		\caption{Likelihood Ratio Test Results}
		\begin{tabular}{cccc}
			\toprule
			\textbf{a} & \textbf{p} & $-2\log\Lambda_1$ & \textbf{Decision on $H_0$} \\ \midrule
			0.1 & 0.1 & 1.9648 & Fail to reject \\ 
			0.1 & 0.2 & 75.3318 & Reject \\ 
			0.1 & 0.3 & 465.3507 & Reject \\ 
			0.1 & 0.4 & 1302.2148 & Reject \\ \midrule
			0.2 & 0.1 & 62.5624 & Reject \\ 
			0.2 & 0.2 & 0.5426 & Fail to reject \\ 
			0.2 & 0.3 & 91.6928 & Reject \\ 
			0.2 & 0.4 & 537.1743 & Reject \\ \midrule
			0.3 & 0.1 & 178.0430 & Reject \\ 
			0.3 & 0.2 & 85.5790 & Reject \\ 
			0.3 & 0.3 & -0.6683 & Fail to reject \\ 
			0.3 & 0.4 & 143.3914 & Reject \\ \midrule
			0.4 & 0.1 & 315.6163 & Reject \\ 
			0.4 & 0.2 & 282.6401 & Reject \\ 
			0.4 & 0.3 & 93.5413 & Reject \\ 
			0.4 & 0.4 & -0.7403 & Fail to reject \\ \midrule
			0.5 & 0.1 & 598.6498 & Reject \\ 
			0.5 & 0.2 & 678.1999 & Reject \\ 
			0.5 & 0.3 & 443.8051 & Reject \\ 
			0.5 & 0.4 & 175.6443 & Reject \\ \midrule
			0.6 & 0.1 & 986.4705 & Reject \\ 
			0.6 & 0.2 & 1060.7808 & Reject \\ 
			0.6 & 0.3 & 908.5222 & Reject \\ 
			0.6 & 0.4 & 520.2437 & Reject \\ \midrule
			0.7 & 0.1 & 1243.0918 & Reject \\ 
			0.7 & 0.2 & 1596.5124 & Reject \\ 
			0.7 & 0.3 & 1581.8311 & Reject \\ 
			0.7 & 0.4 & 1243.7897 & Reject \\ \midrule
			0.8 & 0.1 & 1548.5890 & Reject \\ 
			0.8 & 0.2 & 2291.3934 & Reject \\ 
			0.8 & 0.3 & 2551.6879 & Reject \\ 
			0.8 & 0.4 & 2361.9013 & Reject \\ \midrule
			0.9 & 0.1 & 2466.6171 & Reject \\ 
			0.9 & 0.2 & 3305.7161 & Reject \\ 
			0.9 & 0.3 & 3631.3346 & Reject \\ 
			0.9 & 0.4 & 3835.4538 & Reject \\ 
			\bottomrule
		\end{tabular}
		\label{chi_squared_results1}
	\end{minipage}%
	\hfill
	\begin{minipage}{0.45\textwidth}
		\centering
		\setlength{\tabcolsep}{6pt} % Match column separation with the other table
		\renewcommand{\arraystretch}{1.2} % Adjust row height to match
		\caption{KS Distance values}
		\begin{tabular}{cccc}
			\hline
			$a$ & $p$ & KS Distance &Decision on $H_0$\\ \hline
			0.1 & 0.1 & 0.0027 & Fail to reject\\ 
			0.1 & 0.2 & 0.0181& Reject \\ 
			0.1 & 0.3 & 0.0603 & Reject\\ 
			0.1 & 0.4 & 0.1188 & Reject\\ \midrule
			0.2 & 0.1 & 0.0111 & Reject\\ 
			0.2 & 0.2 & 0.0009 & Fail to reject \\ 
			0.2 & 0.3 & 0.0297 & Reject\\ 
			0.2 & 0.4 & 0.0775 & Reject\\ \midrule
			0.3 & 0.1 & 0.0203 & Reject\\ 
			0.3 & 0.2 & 0.0211 & Reject\\ 
			0.3 & 0.3 & 0.0013 & Fail to reject\\ 
			0.3 & 0.4 & 0.0433 & Reject\\ \midrule
			0.4 & 0.1 & 0.0284 & Reject\\ 
			0.4 & 0.2 & 0.0403 & Reject\\ 
			0.4 & 0.3 & 0.0286 & Reject\\ 
			0.4 & 0.4 & 0.0020 & Fail to reject\\ \midrule
			0.5 & 0.1 & 0.0379 & Reject\\ 
			0.5 & 0.2 & 0.0611 & Reject\\ 
			0.5 & 0.3 & 0.0617 & Reject\\ 
			0.5 & 0.4 & 0.0405 & Reject\\ \midrule
			0.6 & 0.1 & 0.0507 & Reject\\ 
			0.6 & 0.2 & 0.0786 & Reject\\ 
			0.6 & 0.3 & 0.0887 & Reject\\ 
			0.6 & 0.4 & 0.0781 & Reject\\ \midrule
			0.7 & 0.1 & 0.0606 & Reject\\ 
			0.7 & 0.2 & 0.0992 & Reject\\ 
			0.7 & 0.3 & 0.1177 & Reject\\ 
			0.7 & 0.4 & 0.1200 & Reject\\ \midrule
			0.8 & 0.1 & 0.0701 & Reject\\ 
			0.8 & 0.2 & 0.1189 & Reject\\ 
			0.8 & 0.3 & 0.1510 & Reject\\ 
			0.8 & 0.4 & 0.1623 & Reject\\ \midrule
			0.9 & 0.1 & 0.0916 & Reject\\ 
			0.9 & 0.2 & 0.1445 & Reject\\ 
			0.9 & 0.3 & 0.1718 & Reject\\ 
			0.9 & 0.4 & 0.1969 & Reject\\ \bottomrule
		\end{tabular}
		\label{ks1}
	\end{minipage}
\end{table}

\begin{table}[htbp]
	\centering
	%\caption*{Likelihood Ratio Test Results and KS Distance values for different combinations of $a$ and $p$}
	\begin{minipage}{0.5\textwidth}
		\centering
		\setlength{\tabcolsep}{6pt} % Adjust column separation to match both tables
		\renewcommand{\arraystretch}{1.2} % Adjust row height
		\caption{Likelihood Ratio Test Results}
		\begin{tabular}{cccc}
			\toprule
			\textbf{a} & \textbf{p} & 
			$-2\log\Lambda_2$ & \textbf{Decision on $H_0$} \\
			\midrule
			0.1 & 0.6 & 279.8288 & Reject \\
			0.1 & 0.7 & 107.8140 & Reject \\
			0.1 & 0.8 & 16.0485 & Reject \\
			0.1 & 0.9 & 0.3192 & Fail to reject \\ \midrule
			0.2 & 0.6 & 111.3087 & Reject \\
			0.2 & 0.7 & 17.5196 & Reject \\
			0.2 & 0.8 & 0.4450 & Fail to reject \\
			0.2 & 0.9 & 12.1072 & Reject \\ \midrule
			0.3 & 0.6 & 38.3697 & Reject \\
			0.3 & 0.7 & -1.3511 & Fail to reject \\
			0.3 & 0.8 & 15.9764 & Reject \\
			0.3 & 0.9 & 26.5985 & Reject \\ \midrule
			0.4 & 0.6 & 0.0341 & Fail to reject \\
			0.4 & 0.7 & 17.6239 & Reject \\
			0.4 & 0.8 & 50.2535 & Reject \\
			0.4 & 0.9 & 62.1340 & Reject \\ \midrule
			0.5 & 0.6 & 16.9398 & Reject \\
			0.5 & 0.7 & 97.3027 & Reject \\
			0.5 & 0.8 & 127.9626 & Reject \\
			0.5 & 0.9 & 165.1340 & Reject \\ \midrule
			0.6 & 0.6 & 118.6838 & Reject \\
			0.6 & 0.7 & 195.6122 & Reject \\
			0.6 & 0.8 & 209.1521 & Reject \\
			0.6 & 0.9 & 222.3218 & Reject \\ \midrule
			0.7 & 0.6 & 279.7745 & Reject \\
			0.7 & 0.7 & 328.9172 & Reject \\
			0.7 & 0.8 & 315.8419 & Reject \\
			0.7 & 0.9 & 263.6996 & Reject \\ \midrule
			0.8 & 0.6 & 522.4635 & Reject \\
			0.8 & 0.7 & 557.5287 & Reject \\
			0.8 & 0.8 & 533.2392 & Reject \\
			0.8 & 0.9 & 284.6627 & Reject \\ \midrule
			0.9 & 0.6 & 768.7074 & Reject \\
			0.9 & 0.7 & 688.8740 & Reject \\
			0.9 & 0.8 & 545.8894 & Reject \\
			0.9 & 0.9 & 411.5155 & Reject \\   
			\bottomrule
		\end{tabular}
		\label{chi_squared_results2}
	\end{minipage}%
	\hfill
	\begin{minipage}{0.45\textwidth}
		\centering
		\setlength{\tabcolsep}{6pt} % Match column separation with the other table
		\renewcommand{\arraystretch}{1.2} % Adjust row height to match
		\caption{KS Distance values}
		\begin{tabular}{cccc}
			\hline
			$a$ & $p$ & KS Distance&Decision on $H_0$ \\ \hline
			0.1 & 0.6 & 0.1188 & Reject\\
			0.1 & 0.7 & 0.0603 & Reject\\
			0.1 & 0.8 & 0.0181 & Reject\\
			0.1 & 0.9 & 0.0027 & Fail to reject\\ \midrule
			0.2 & 0.6 & 0.0775 & Reject\\
			0.2 & 0.7 & 0.0297 & Reject\\
			0.2 & 0.8 & 0.0009 & Fail to reject\\
			0.2 & 0.9 & 0.0111 & Reject\\
			\midrule
			0.3 & 0.6 & 0.0433 & Reject\\
			0.3 & 0.7 & 0.0013 & Fail to reject\\
			0.3 & 0.8 & 0.0211 & Reject\\
			0.3 & 0.9 & 0.0203 & Reject\\
			\midrule
			0.4 & 0.6 & 0.0020 & Fail to reject\\
			0.4 & 0.7 & 0.0286 & Reject\\
			0.4 & 0.8 & 0.0403 & Reject\\
			0.4 & 0.9 & 0.0284 & Reject\\
			\midrule
			0.5 & 0.6 & 0.0365 & Reject\\
			0.5 & 0.7 & 0.0574 & Reject\\
			0.5 & 0.8 & 0.0559 & Reject\\
			0.5 & 0.9 & 0.0409 & Reject\\
			\midrule
			0.6 & 0.6 & 0.0781 & Reject\\
			0.6 & 0.7 & 0.0887 & Reject\\
			0.6 & 0.8 & 0.0786 & Reject\\
			0.6 & 0.9 & 0.0507 & Reject\\
			\midrule
			0.7 & 0.6 & 0.1200 & Reject\\
			0.7 & 0.7 & 0.1177 & Reject\\
			0.7 & 0.8 & 0.0992 & Reject\\
			0.7 & 0.9 & 0.0606 & Reject\\
			\midrule
			0.8 & 0.6 & 0.1623 & Reject\\
			0.8 & 0.7 & 0.1510 & Reject\\
			0.8 & 0.8 & 0.1189 & Reject\\
			0.8 & 0.9 & 0.0701 & Reject\\
			\midrule
			0.9 & 0.6 & 0.1969 & Reject\\
			0.9 & 0.7 & 0.1718 & Reject\\
			0.9 & 0.8 & 0.1445 & Reject\\
			0.9 & 0.9 & 0.0916 & Reject\\    
			
			\bottomrule
		\end{tabular}
		\label{ks2}
	\end{minipage}
\end{table}

The test indicates that when $a=p<1/2$ or $a=1-p\leq 1/2$, the Markov chain model generated by copula (\ref{our model}) and the Bernoulli($p$) marginal distribution transforms into independent Bernoulli trials.

\subsubsection{Kolmogorov-Smirnov (KS) Distance}
We employed the Kolmogorov-Smirnov (KS) distance to compare the copula of the sample with the independence copula. Specifically, we used it to test whether the empirical joint distribution of $(X_t, X_{t+1})$ deviates significantly from the product of their marginal distributions. The KS distance is given by

\begin{align*}
	KS = \sup_{x} \left| \hat{F}(X_t, X_{t+1}) - \hat{F}(X_t) \cdot \hat{F}(X_{t+1}) \right|,
\end{align*}
where $\hat{F}(X_t, X_{t+1})$ is the empirical joint distribution function of consecutive states, and $\hat{F}(X_t)$ and $\hat{F}(X_{t+1})$ are the empirical marginal distribution functions of $X_t$ and $X_{t+1}$, respectively.

A large KS distance suggests a significant departure from independence, indicating the presence of dependence in the Markov chain. The results of this test are summarized in Tables \ref{ks1} and \ref{ks2}.

For this test, we defined the null hypothesis $H_0$: The consecutive states in the Markov chain are independent, against the alternative hypothesis $H_1$: The consecutive states exhibit some form of dependence. We conducted the test at a significance level of $\alpha = 0.05$ for a sequence of length $n=9999$. The test shows that when $a = p < \frac{1}{2}$ or $a = 1 - p \leq \frac{1}{2}$, the Markov chain model generated by the copula (\ref{our model}) with Bernoulli($p$) marginals reduces to a sequence of independent Bernoulli trials.

\subsection{Comparison of different estimators for the mean} 
This section is dedicated to the performance of the different estimators of the Bernoulli parameter. The MLEs are obtained numerically and their confidence intervals are constructed.  The sample mean satisfies the central limit theorem (\ref{clt of xbar}). The $(1-\alpha)100\%$ confidence intervals for $p$ are therefore \begin{align*}
	\begin{cases}
		\Bar{p}\pm z_{\alpha/2}\left[\frac{\Bar{p}(1-\Bar{p})}{n+1}\left(\frac{1+a-2\Bar{p}}{1-a}\right)\right]^{1/2},& \quad \text{when}\ \bar{p}<1/2\\
		\Bar{p}\pm z_{\alpha/2}\left[\frac{\Bar{p}(1-\Bar{p})}{n+1}\left(\frac{a+2\Bar{p}-1}{1-a}\right)\right]^{1/2}, & \quad \text{when}\ \bar{p}\geq 1/2.     
	\end{cases}
\end{align*}

The other estimator for the mean ($p$) considered here is the robust estimator proposed by Longla and Peligrad (2021)  \cite{longla2021new}. It is based on a sample of dependent observations and requires some mild conditions to be satisfied on the variance of partial sums. The procedure is as follows: Given the sample data $(X_t)_{0\leq t\leq n}$ from a stationary sequence $(X_t)_{t\in\mathbb{Z}}$, we generate another random sample $(Y_t)_{0\leq t\leq n}$ independently of $(X_t)_{0\leq t\leq n}$ and following a standard normal distribution. The Gaussian kernel and the optimal bandwidth $h_n$ are used. The optimal bandwidth, as proposed in Longla and Peligrad (2021)  \cite{longla2021new}, for our case is 
\begin{align}\label{optimal bandwidth}
	h_n=\left(\frac{1}{(n+1)\sqrt{2}}\right)^{1/5}. 
\end{align}

To use this estimator, the following conditions must be satisfied
\textit{(i.)} $(X_t)_{t\in\mathbb{Z}}$ is an ergodic sequence,
\textit{(ii.)} $(X_t)_{t\in\mathbb{Z}}$ has finite second moments and    \textit{(iii.)} $(n+1)h_n var(\Bar{p})\longrightarrow 0$ as $n\longrightarrow\infty.$
These conditions are easy to verify: ($i.$) The Markov chain $X$ is $\psi$-mixing, and hence ergodic. ($ii.$) The second condition follows from that $Var(X_0)=p(1-p)<\infty\ \text{and}\ E[X_0]=p$.  ($iii.$) Using the variance of $\Bar{X}$ given by (\ref{varbarX1}), we see that \begin{align*}
	(n+1)\left(\frac{1}{(n+1)\sqrt{2}}\right)^{1/5}\left[\frac{p(1-p)}{(n+1)}\left(\frac{1+a-2p}{1-a}\right)+O(1/n)\right]\longrightarrow 0
\end{align*}
as \ $n\longrightarrow\infty.$ The proposed robust estimator for the mean $(p)$ of $X$ is 
\begin{equation*}
	\tilde{p}=\frac{1}{(n+1)h_n}\sum_{t=0}^n X_{t}\exp\left(-\frac{1}{2}\left(\frac{Y_t}{h_n}\right)^2\right).
\end{equation*}	
A $(1-\alpha)100\%$ confidence interval for $p$ is\begin{equation*}
	\tilde{p}\sqrt{1+h_n^2}\pm z_{\alpha/2}\left[\frac{\Bar{X}_n}{(n+1)\sqrt{2}h_n}\right]^{1/2}.
\end{equation*}

Tables \ref{table for p<1/2} and \ref{table for p>=1/2} below present a comparison of the three estimators of $p$ across varying true values of $p$. Samples of size $n+1$ were utilized, and 95\% confidence intervals were constructed. Coverage probabilities were computed based on $400$ replications of each of the chains of size $n+1$. Subsequently, the means of the lengths of the $400$ confidence intervals were calculated to determine the confidence interval mean length.  The results in both tables demonstrate that the estimators are asymptotically equivalent and converge to the true value. Note that we did not include the numerical values of the estimators because $400$ numerical values were computed for each estimator in each case, and there is no criterion to choose only one representative from the $400$. Relatively narrow confidence intervals reflect a higher level of precision in estimation, while wider confidence intervals indicate greater uncertainty in estimation.

%\bgroup
%\begingroup

\setlength{\tabcolsep}{6pt} % Default value: 6pt
\def\arraystretch{1.0}
\begin{table}[h]
	\caption{A comparison of the performance of  $\hat{p}$, $\Bar{p}$ and  $\Tilde{p}$ for $p<0.5$. A fixed value of $a=0.5$ was used in all the cases and 95\% confidence intervals were constructed. \ \textit{Est} denotes estimator and $p_0$ denotes the true value of $p$.}
	\label{table for p<1/2}
	\centering
	{\scriptsize     \begin{tabular}{|c||*{5}{c|}}\hline
			n & \backslashbox{Est}{$p_0$} &\makebox[2em]{0.1}&\makebox[2em]{0.2}&\makebox[2em]{0.3} & 0.4\\ \hline\hline
			
			\multirow{2}{*}{99}    & $\hat{p}$ & \begin{tabular}{c}
				CIML: 0.1862\\
				CP: 0.8625
			\end{tabular} & \begin{tabular}{c}
				CIML: 0.2258\\
				CP: 0.9100
			\end{tabular} &\begin{tabular}{c}
				CIML: 0.2355\\
				CP: 0.9275
			\end{tabular} & \begin{tabular}{c}
				CIML: 0.2245\\
				CP: 0.9400
			\end{tabular} \\
			\cline{2-6}
			& $\Bar{p}$ & \begin{tabular}{c}
				CIML: 0.1896\\
				CP: 0.9625
			\end{tabular} & \begin{tabular}{c}
				CIML:0.0.2325\\
				CP: 0.9475
			\end{tabular} & \begin{tabular}{c}
				CIML: 0.2410\\
				CP: 0.9575
			\end{tabular}& \begin{tabular}{c}
				CIML: 0.2272\\
				CP: 0.9625
			\end{tabular} \\
			\cline{2-6}
			& $\Tilde{p}$  &\begin{tabular}{c}
				CIML: 0.1656\\
				CP: 0.7875
			\end{tabular}  & \begin{tabular}{c}
				CIML: 0.2399\\
				CP: 0.9000
			\end{tabular} &\begin{tabular}{c}
				CIML: 0.2946\\
				CP: 0.9325
			\end{tabular} &\begin{tabular}{c}
				CIML: 0.3416\\
				CP: 0.9500
			\end{tabular}  \\
			\hline
			\hline
			\multirow{2}{*}{499}    & $\hat{p}$ & \begin{tabular}{c}
				CIML: 0.0828\\
				CP: 0.9050
			\end{tabular} & \begin{tabular}{c}
				CIML: 0.1030\\
				CP: 0.9150
			\end{tabular} &\begin{tabular}{c}
				CIML: 0.1074\\
				CP: 0.9450
			\end{tabular} & \begin{tabular}{c}
				CIML:0.1012\\
				CP: 0.9300
			\end{tabular} \\
			\cline{2-6}
			& $\Bar{p}$ & \begin{tabular}{c}
				CIML: 0.0848\\
				CP: 0.9450
			\end{tabular} & \begin{tabular}{c}
				CIML: 0.1040\\
				CP: 0.9550
			\end{tabular} &\begin{tabular}{c}
				CIML: 0.1078\\
				CP: 0.9375
			\end{tabular} &\begin{tabular}{c}
				CIML:0.1016\\
				CP: 0.9275
			\end{tabular}  \\
			\cline{2-6}
			& $\Tilde{p}$  &\begin{tabular}{c}
				CIML: 0.0889\\
				CP: 0.8825
			\end{tabular}  & \begin{tabular}{c}
				CIML: 0.1264\\
				CP: 0.9325
			\end{tabular} &\begin{tabular}{c}
				CIML: 0.1550\\
				CP: 0.9250
			\end{tabular} &\begin{tabular}{c}
				CIML: 0.1797\\
				CP: 0.9675
			\end{tabular} \\
			\hline
			\hline
			\multirow{2}{*}{999}    & $\hat{p}$ & \begin{tabular}{c}
				CIML: 0.0591\\
				CP: 0.9375
			\end{tabular} & \begin{tabular}{c}
				CIML: 0.0732\\
				CP: 0.9400
			\end{tabular} & \begin{tabular}{c}
				CIML: 0.0761\\
				CP: 0.9375
			\end{tabular}& \begin{tabular}{c}
				CIML: 0.0717\\
				CP: 0.9475
			\end{tabular} \\
			\cline{2-6}
			& $\Bar{p}$ & \begin{tabular}{c}
				CIML: 0.0600\\
				CP: 0.9375
			\end{tabular} & \begin{tabular}{c}
				CIML: 0.0735\\
				CP: 0.9700
			\end{tabular} & \begin{tabular}{c}
				CIML: 0.0762\\
				CP: 0.9575
			\end{tabular}&  \begin{tabular}{c}
				CIML: 0.0719\\
				CP: 0.9550
			\end{tabular}\\
			\cline{2-6}
			& $\Tilde{p}$  &\begin{tabular}{c}
				CIML: 0.0678\\
				CP: 0.8825
			\end{tabular}  & \begin{tabular}{c}
				CIML: 0.0962\\
				CP: 0.9200
			\end{tabular} &\begin{tabular}{c}
				CIML: 0.1179\\
				CP: 0.9300
			\end{tabular} &\begin{tabular}{c}
				CIML: 0.1362\\
				CP: 0.9550
			\end{tabular}   \\
			\hline
			\hline
			\multirow{2}{*}{4999}    & $\hat{p}$ & \begin{tabular}{c}
				CIML: 0.0267\\
				CP: 0.9425
			\end{tabular} & \begin{tabular}{c}
				CIML: 0.0328\\
				CP: 0.9575
			\end{tabular} &\begin{tabular}{c}
				CIML: 0.0341\\
				CP: 0.9475
			\end{tabular} & \begin{tabular}{c}
				CIML: 0.0321\\
				CP: 0.9350
			\end{tabular} \\
			\cline{2-6}
			& $\Bar{p}$ & \begin{tabular}{c}
				CIML: 0.0268\\
				CP: 0.9625
			\end{tabular} &\begin{tabular}{c}
				CIML: 0.0329\\
				CP: 0.9600
			\end{tabular}  &\begin{tabular}{c}
				CIML: 0.0341\\
				CP: 0.9400
			\end{tabular} & \begin{tabular}{c}
				CIML: 0.0321\\
				CP: 0.9525
			\end{tabular} \\
			\cline{2-6}
			& $\Tilde{p}$  &\begin{tabular}{c}
				CIML:0.0356\\
				CP: 0.9000
			\end{tabular}  & \begin{tabular}{c}
				CIML: 0.0505\\
				CP: 0.9200
			\end{tabular} &\begin{tabular}{c}
				CIML: 0.0619 \\
				CP: 0.9425
			\end{tabular} &\begin{tabular}{c}
				CIML: 0.0716\\
				CP: 0.9450
			\end{tabular}  \\
			\hline
			\multirow{2}{*}{9999}    & $\hat{p}$ & \begin{tabular}{c}
				CIML: 0.0189\\
				CP: 0.9525
			\end{tabular} & \begin{tabular}{c}
				CIML: 0.0232\\
				CP: 0.9550
			\end{tabular} &\begin{tabular}{c}
				CIML: 0.0241\\
				CP: 0.9500
			\end{tabular} & \begin{tabular}{c}
				CIML: 0.0227\\
				CP: 0.9300
			\end{tabular} \\
			\cline{2-6}
			& $\Bar{p}$ & \begin{tabular}{c}
				CIML: 0.0190 \\
				CP:  0.9500
			\end{tabular} & \begin{tabular}{c}
				CIML: 0.0232\\
				CP: 0.9570
			\end{tabular} & \begin{tabular}{c}
				CIML: 0.0241\\
				CP: 0.9625
			\end{tabular}&  \begin{tabular}{c}
				CIML: 0.0227\\
				CP: 0.9450
			\end{tabular}\\
			\cline{2-6}
			& $\Tilde{p}$  &\begin{tabular}{c}
				CIML: 0.0270 \\
				CP: 0.9025
			\end{tabular}  & \begin{tabular}{c}
				CIML: 0.0383\\
				CP: 0.9000
			\end{tabular} &\begin{tabular}{c}
				CIML: 0.0470\\
				CP: 0.9500
			\end{tabular} &\begin{tabular}{c}
				CIML: 0.0542\\
				CP: 0.9500
			\end{tabular}  \\
			\hline
			\hline    
	\end{tabular} }
\end{table}
%\egroup

 \bgroup
\def\arraystretch{1.0}
\begin{table}[htbp]
	\caption{\small A comparison of the performance of  $\hat{p}$, $\Bar{p}$ and  $\Tilde{p}$ for $p\geq 0.5$. A fixed value $a=0.5$ was used in all the cases and 95\% confidence intervals were constructed. \ \textit{Est} denotes estimator and $p_0$ denotes the true value of $p$.}
	\label{table for p>=1/2}
	\centering
	{\scriptsize \begin{tabular}{|l||*{5}{c|}}\hline
			n & \backslashbox{Est}{$p_0$}
			& \makebox[1.5em]{0.6}&\makebox[1.5em]{0.7}&\makebox[1.5em]{0.8}&\makebox[1.5em]{0.9}\\ \hline\hline
			
			\multirow{2}{*}{99}    & $\hat{p}$ & \begin{tabular}{c}
				CIML: 0.2242\\
				CP: 0.9575
			\end{tabular} &\begin{tabular}{c}
				CIML: 0.2381\\
				CP: 0.9275
			\end{tabular} & \begin{tabular}{c}
				CIML: 0.2286                                  \\
				CP: 0.9250
			\end{tabular} &\begin{tabular}{c}
				
				CIML:    0.1912                              \\
				CP: 0.8900
			\end{tabular}\\
			\cline{2-6}
			& $\Bar{p}$ & \begin{tabular}{c}
				CIML: 0.2272\\
				CP: 0.9450
			\end{tabular} &\begin{tabular}{c}
				CIML: 0.2410\\
				CP: 0.9575
			\end{tabular} & \begin{tabular}{c}
				CIML: 0.2326\\
				CP: 0.9450
			\end{tabular} & \begin{tabular}{c}
				CIML: 0.1896\\
				CP: 0.9500
			\end{tabular}\\
			\cline{2-6}
			& $\Tilde{p}$  & \begin{tabular}{c}
				CIML: 0.4163\\
				CP: 0.9450
			\end{tabular} &\begin{tabular}{c}
				CIML: 0.4529\\
				CP: 0.9525
			\end{tabular} &\begin{tabular}{c}
				CIML: 0.4841\\
				CP: 0.9575
			\end{tabular} & \begin{tabular}{c}
				CIML: 0.5124\\
				CP: 0.9700
			\end{tabular} \\
			\hline
			\hline
			\multirow{2}{*}{499}    & $\hat{p}$ &  \begin{tabular}{c}
				CIML: 0.1013\\
				CP: 0.9200
			\end{tabular} &\begin{tabular}{c}
				
				CIML: 0.1071\\
				CP: 0.9425
			\end{tabular} & \begin{tabular}{c}
				CIML: 0.1032                             \\
				CP: 0.9550
			\end{tabular} &\begin{tabular}{c}
				
				CIML:  0.0832                                      \\
				CP: 0.9075
			\end{tabular} \\
			\cline{2-6}
			& $\Bar{p}$ & \begin{tabular}{c}
				CIML: 0.1016\\
				CP: 0.9575
			\end{tabular} &\begin{tabular}{c}
				CIML: 0.1077\\
				CP: 0.9475
			\end{tabular} &\begin{tabular}{c}
				CIML:  0.1040\\
				CP: 0.9475
			\end{tabular} &\begin{tabular}{c}
				CIML: 0.0848\\
				CP: 0.9400
			\end{tabular} \\
			\cline{2-6}
			& $\Tilde{p}$  & \begin{tabular}{c}
				
				CIML: 0.2200\\
				CP: 0.9450
			\end{tabular} &\begin{tabular}{c}
				CIML: 0.2381\\
				CP: 0.9650
			\end{tabular} &\begin{tabular}{c}
				CIML: 0.2543\\
				CP: 0.9500
			\end{tabular} & \begin{tabular}{c}
				CIML: 0.2695\\
				CP: 0.9775
			\end{tabular}\\
			\hline
			\hline
			\multirow{2}{*}{999}    & $\hat{p}$  & \begin{tabular}{c}
				CIML: 0.0716\\
				CP: 0.9275
			\end{tabular} & \begin{tabular}{c}
				CIML: 0.0759\\
				CP: 0.9175
			\end{tabular}& \begin{tabular}{c}
				CIML:    0.0731                                   \\
				CP: 0.9400
			\end{tabular} &\begin{tabular}{c}
				CIML:  0.0595                                    \\
				CP: 0.9225
			\end{tabular}\\
			\cline{2-6}
			&  $\Bar{p}$&  \begin{tabular}{c}
				CIML: 0.0719\\
				CP: 0.9475
			\end{tabular} &\begin{tabular}{c}
				CIML: 0.0762\\
				CP: 0.9350
			\end{tabular} &\begin{tabular}{c}
				CIML: 0.0735\\
				CP: 0.9625
			\end{tabular} &\begin{tabular}{c}
				CIML: 0.0600\\
				CP: 0.9175
			\end{tabular} \\
			\cline{2-6}
			& $\Tilde{p}$ & \begin{tabular}{c}
				CIML: 0.1664\\
				CP: 0.9525
			\end{tabular} &\begin{tabular}{c}
				CIML: 0.1803\\
				CP: 0.9725
			\end{tabular} &\begin{tabular}{c}
				CIML: 0.1926\\
				CP: 0.9550
			\end{tabular} & \begin{tabular}{c}
				CIML: 0.2042\\
				CP: 0.9375
			\end{tabular} \\
			\hline
			\hline
			\multirow{2}{*}{4999}    & $\hat{p}$  &\begin{tabular}{c}
				CIML: 0.0321\\
				CP: 0.9400
			\end{tabular} & \begin{tabular}{c}
				CIML: 0.0341\\
				CP: 0.9425
			\end{tabular} &\begin{tabular}{c}
				CIML:  0.0329                                       \\
				CP: 0.9500
			\end{tabular} & \begin{tabular}{c}
				CIML:  0.0268                                       \\
				CP: 0.9400
			\end{tabular} \\
			\cline{2-6}
			& $\Bar{p}$  & \begin{tabular}{c}
				CIML: 0.0321\\
				CP: 0.9550
			\end{tabular} & \begin{tabular}{c}
				CIML: 0.0341\\
				CP: 0.9275
			\end{tabular}& \begin{tabular}{c}
				CIML: 0.0329\\
				CP: 0.9525
			\end{tabular}& \begin{tabular}{c}
				CIML: 0.0268\\
				CP: 0.9625
			\end{tabular}\\
			\cline{2-6}
			& $\Tilde{p}$  & \begin{tabular}{c}
				CIML: 0.0876\\
				CP: 0.9475
			\end{tabular} &\begin{tabular}{c}
				CIML: 0.0946\\
				CP: 0.9525
			\end{tabular} &\begin{tabular}{c}
				CIML: 0.1012\\
				CP: 0.9500
			\end{tabular} & \begin{tabular}{c}
				CIML: 0.1073\\
				CP: 0.9525
			\end{tabular} \\
			\hline
			\hline
			\multirow{2}{*}{9999}    & $\hat{p}$ & \begin{tabular}{c}
				CIML: 0.0227 \\
				CP:  0.9500
			\end{tabular} &\begin{tabular}{c}
				CIML: 0.0241 \\
				CP:  0.9500
			\end{tabular} & \begin{tabular}{c}
				CIML: 0.0233 \\
				CP:  0.9475
			\end{tabular} &\begin{tabular}{c}
				CIML: 0.0190 \\
				CP:  0.9425
			\end{tabular}\\
			\cline{2-6}
			& $\Bar{p}$ &  \begin{tabular}{c}
				CIML: 0.0227 \\
				CP:  0.9500
			\end{tabular}& \begin{tabular}{c}
				CIML:  0.0241\\
				CP:  0.9450
			\end{tabular}& \begin{tabular}{c}
				CIML: 0.0233 \\
				CP:  0.9450
			\end{tabular}& \begin{tabular}{c}
				CIML: 0.0190 \\
				CP:  0.9500
			\end{tabular}\\
			\cline{2-6}
			& $\Tilde{p}$  & \begin{tabular}{c}
				CIML: 0.0664 \\
				CP:  0.9700
			\end{tabular} &\begin{tabular}{c}
				CIML: 0.0717 \\
				CP:  0.9725
			\end{tabular} &\begin{tabular}{c}
				CIML: 0.0767 \\
				CP:  0.9700
			\end{tabular} & \begin{tabular}{c}
				CIML:  0.0813\\
				CP:  0.9650
			\end{tabular} \\
			\hline
			\hline    
	\end{tabular} }
\end{table}
\egroup

\begin{figure}[htbp]
    \centering

    \begin{subfigure}{0.5\textwidth}
        \centering
        \fbox{\includegraphics[width=\textwidth]{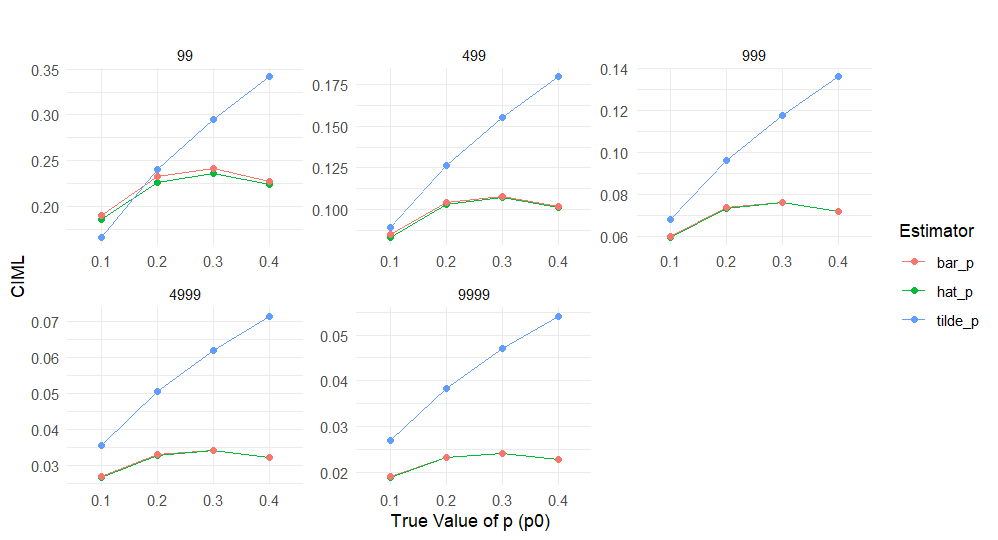}}
        \caption{CIML for estimators when $p < \frac{1}{2}$.}
        \label{CIML-less}
    \end{subfigure}

    \vspace{1em}

    \begin{subfigure}{0.5\textwidth}
        \centering
        \fbox{\includegraphics[width=\textwidth]{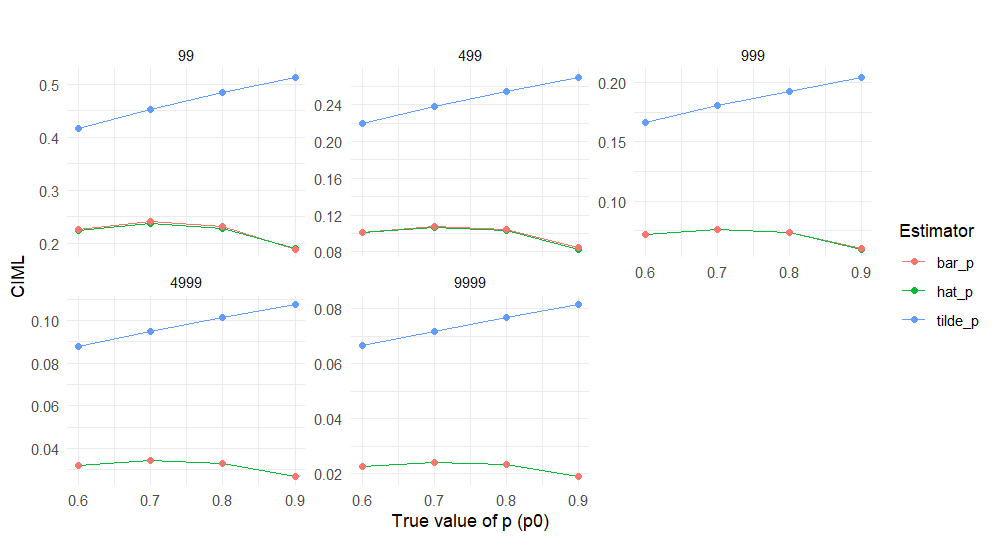}}
        \caption{CIML for estimators when $p >\frac{1}{2}$.}
        \label{CIML-geq}
    \end{subfigure}

    \caption{Confidence Interval Mean Lengths (CIML) for different values of $p$.}
    \label{CIML}
\end{figure}

In both cases of $p<1/2$ and $p> 1/2$, $\hat{p}$ yields the most precise estimates with the shortest mean confidence intervals. The sample mean closely follows $\hat{p}$ in precision but exhibits slightly longer mean confidence intervals. Although the robust estimator shows larger mean confidence intervals compared to $\hat{p}$ and $\bar{p}$, it still provides reasonable estimates, especially with larger sample sizes. This comparison is presented on the plots in Figure \ref{CIML-less} and Figure \ref{CIML-geq}.

Both the MLE and the robust estimator demonstrate sensitivity to sample size, resulting in initially low coverage probabilities. However, as the sample size increases, these probabilities approach the desired 95\% level. The coverage probabilities for the sample mean $(\bar{p})$ consistently approach the desired nominal level, indicating reliable confidence intervals across different sample sizes.

The performance of the robust estimator is less optimal with smaller sample sizes but improves as the sample size increases. The robust estimator is designed to be competitive in data containing noise or outliers; however, the data considered here has neither.

While all three estimators provide asymptotically unbiased estimates for the Bernoulli parameter $p$, the sample mean ($\bar{p}$) generally strikes the best balance between precision and accuracy, making it the preferred choice for estimating $p$ in this context.

\clearpage      % flush all floats (tables/figures)
\section{Conclusion}
The mixing properties of copula-based Markov chains are highly dependent on the chosen marginal distributions. Some copulas that do not generate mixing with continuous marginals may generate mixing with discrete marginals, so generalizations should be approached with caution. In this paper, we have demonstrated that copulas from the Fr\'echet (Mardia) family generate $\psi$-mixing Markov chains for $a+b=1$, which is not the case when the marginals are uniform$[0,1]$. Understanding these mixing properties is crucial for deriving central limit theorems for sample means, which in turn facilitates the construction of confidence intervals for the sample means based on the standard normal distribution.

\begin{funding}
 Authors state no funding involved.  
\end{funding}

\begin{authorcontributions}
	All authors share full responsibility for the content of this manuscript. They have jointly agreed to its submission, carefully examined the results, and approved the final version for publication.
\end{authorcontributions}

\begin{conflictofinterest}
	The authors declare that they have no competing interests.
\end{conflictofinterest}

\end{document}